\documentclass{amsart}
\usepackage{amssymb,latexsym}
\usepackage{amscd,amsthm}

\usepackage[all]{xy}
\usepackage{tikz}

\usepackage{tikz-cd}

\newenvironment{psmallmatrix}
{\left(\begin{smallmatrix}}
	{\end{smallmatrix}\right)}

\newtheorem{thm}{\bf Theorem}[section]
\newtheorem{cor}[thm]{\bf Corollary}

\newtheorem{lem}[thm]{\bf Lemma}
\newtheorem{prop}[thm]{\bf Proposition}
\newtheorem{defn}[thm]{\bf Definition}

\newtheorem{rem}[thm]{\bf Remark}

\newtheorem{exmp}[thm]{\bf Example}

\DeclareMathOperator{\Ext}{Ext}
\DeclareMathOperator{\Hom}{Hom}

\newcommand{\N}{\mathbb{N}}

\newcommand{\cqfd}
{\hspace{1cm}
\rule{2mm}{2mm}%
\medbreak%
\par%
}

\def\pd{{\rm pd}}

\def\id{{\rm id}}

\def\Mor{{\rm Mor}}

\def\GP{\mathcal{GP}}
\def\GI{\mathcal{GI}}

\def\A{\mathcal{A}}
\def\B{\mathcal{B}}
\def\C{\mathcal{C}}

\def\T{\mathcal{T}}

\def\X{\mathcal{X}}
\def\Y{\mathcal{Y}}

\def\PP{\rm \textbf{P}}
\def\EE{\rm \textbf{E}}

\def\F{\mathcal{F}}
\def\P{\mathcal{P}}
\def\I{\mathcal{I}}

\def\h{{\rm \textbf{h}}}

\def\pp{{\rm \textbf{p}}}
\def\hh{{\rm \textbf{h}}}
\def\qq{{\rm \textbf{q}}}
\def\zz{{\rm \textbf{z}}}
\def\cc{{\rm \textbf{c}}}
\def\kk{{\rm \textbf{k}}}

\def\PP{{\rm \textbf{P}}}
\def\QQ{{\rm \textbf{Q}}}

\def\G{{\rm G}}

\def\Mor{{\rm Mor}}

\def\resdim{{\rm resdim}}
\def\coresdim{{\rm coresdim}}

\def\glresdim{{\rm gl.resdim}}
\def\glcoresdim{{\rm gl.coresdim}}

\def\Gpd{{\rm Gpd}}

\def\Gid{{\rm Gid}}

\def\glGdim{{\rm gl.Gdim}}

\def\PD{{\rm  PD}}
\def\ID{{\rm  ID}}

\def\lGPD{{\rm  lGPD}}
\def\lGID{{\rm  lGID}}

\def\GPD{{\rm  GPD}}
\def\GID{{\rm  GID}}

\def\Coker{{\rm Coker}}
\def\Ker{{\rm Ker}}

\def\Ext{{\rm Ext}}

\def\Hom{{\rm Hom}}

\begin{document}

\title{Gorenstein homological algebra over comma categories}

\thanks{$^*$ Corresponding author: elmaaouy.rachid@gmail.com}
\author{Rachid El Maaouy$^*$}
\address{SMARTiLab Laboratory, Moroccan School of Engineering Sciences, EMSI Rabat, Morocco}

\email[Rachid El Maaouy]{elmaaouy.rachid@gmail.com}

\author{J. R. García Rozas}

\address{Departamento de Matem\'{a}ticas, Universidad de Almer\'{\i}a, 04071 Almer\'{\i}a, Spain}

\email[J. R. García Rozas]{jrgrozas@ual.es}

\author{Luis Oyonarte}

\address{Departamento de Matem\'{a}ticas, Universidad de Almer\'{\i}a, 04071 Almer\'{\i}a, Spain}
\email[Luis Oyonarte]{oyonarte@ual.es}


\keywords{totally compatible functor, totally acyclic complex, Gorenstein projective dimension, triangular matrix ring, morphism category}

\subjclass[2020]{18A25, 18G25.}

\begin{abstract}

Let $\A$ and $\B$ be abelian categories with enough projective objects. In this paper, we characterize those functors $F:\A\to\B$ for which the Gorenstein projective objects in the comma category $(F\downarrow\B)$ admit an explicit description in terms of the Gorenstein projective objects of $\A$ and $\B$. As a main application, we obtain estimates for the global Gorenstein projective dimension of $(F\downarrow\B)$ in terms of those of $\A$, $\B$ and the functor $F$, and derive exact formulas in several important cases. Our categorical approach also yields dual results concerning Gorenstein injective objects and dimension. Special attention is devoted to the  module category over a triangular matrix ring $
\begin{psmallmatrix}
A & 0\\
U & B
\end{psmallmatrix}$, 
the product category $\A\times\B$, and the $n$-morphism category $\Mor_n(\C)$ associated with an abelian category $\C$ and an integer $n\geq 2$.

\end{abstract}

\maketitle

\section{Introduction}

Triangular matrix rings appear naturally in many areas of mathematics. In particular, they play an important role in ring theory and the representation theory of algebras. They are generally noncommutative and lack left-right symmetry, so they are frequently used to establish examples and counterexamples that enrich the general theories of rings and modules; see, e.g., \cite{FGR75,ARS95}. 
Our present work is motivated by a series of results on triangular matrix rings and, more generally, comma categories. Before discussing in detail the main motivations behind this work, let us introduce some necessary notation. 

Let $A$ and $B$ be rings, and let $U$ be a $(B,A)$-bimodule. There is a ring   $$T=\begin{psmallmatrix}
      A & 0\\ 
      U & B  
    \end{psmallmatrix}=\{\begin{psmallmatrix}
      a & 0\\ 
      u & b  
    \end{psmallmatrix} \mid a\in A, b\in B, u\in U \}$$
called the triangular matrix ring with the natural addition and multiplication. Recall also that a (left) module over $T$ is a triple $M=\left(\begin{smallmatrix}
	 M_1 \\ M_2  \end{smallmatrix}\right) _{\varphi^M}$, where $M_1$ is an $A$-module, $M_2$ is a $B$-module and $\varphi^M : U \otimes_A M_1\rightarrow M_2$ is a $B$-morphism. When studying dual properties over $T$, the morphism $\varphi^M$ can equivalently be replaced by the adjoint morphism $\psi^M: M_1\to \Hom_B(U,M_2)$. A fundamental and natural problem is to determine precise relationships between the homological properties of a $T$-module $M$ and those of $M_1$, $M_2$, and the connecting morphism $\varphi^M$ (equivalently, $\psi^M$). In particular, one seeks relationships between the homological properties of the ring $T$ and those of $A$, $B$, and the bimodule $U$. Classical modules, such as projective and injective modules over $T$, are well described: It follows from the work of   Haghany and Varadarajan \cite{HV99,HV00} that the classes of projective and injective modules over $T$ admit the descriptions  $\P(T)=\Phi(\P(A),\P(B))$ and $\I(T)=\Psi(\I(A),\I(B))$, where for classes $\X_1,\Y_1$ of $A$-modules and classes $\X_2,\Y_2$ of $B$-modules, the class $\Phi(\X_1,\X_2)$ is defined to be the class of all $T$-modules $M=\begin{psmallmatrix}
    M_1\\ M_2
\end{psmallmatrix}_{\varphi^M}$ such that $M_1\in \X_1,\, \Coker\varphi^M\in \X_2$ and $\varphi^M$ is a monomorphism, while $\Psi(\Y_1,\Y_2)$ is defined to be the class of all $T$-modules $M=\begin{psmallmatrix}
    M_1\\ M_2
\end{psmallmatrix}_{\psi^M}$ such that $\Ker \psi^M\in \Y_1, M_2\in \Y_2$ and  $\psi^M$ is an epimorphism.

\medskip

As shown by Zhang in \cite[Examples 1 and 2]{Zha13}, these characterizations do not, in general, extend to the realm of Gorenstein homological algebra. The problem of describing Gorenstein projective and injective modules over triangular matrix rings and algebras has been investigated by many authors under different assumptions (see \cite{Zha13,ECT14,EHSL16,LZHZ20,PZH22,Mao22,HZ22,YHW24}). However, all the assumptions imposed in these works are sufficient but not necessary. This leads us to the following natural question:

\bigskip

\textbf{Question 1:} Given two rings $A$ and $B$, under what assumptions on the $(B,A)$-bimodule $U$ do we have the following characterizations?
$$\GP(T)=\Phi(\GP(A),\GP(B)), \quad \GI(T)=\Psi(\GI(A),\GI(B)).$$
Here, $\GP(R)$ and $\GI(R)$ denote, respectively, the classes of Gorenstein projective and injective modules over a ring $R$.

A partial answer to this question has been given in \cite{Zha13} for an artin algebra $T$. Zhang in \cite{Zha13} introduced compatible bimodules and showed that if $U$ is compatible, then the Gorenstein projective modules over $T$ are explicitly described as above. Moreover, he showed that compatibility is also necessary whenever $A$ is a Gorenstein algebra. Our first main goal of this paper is to develop this idea further, giving a positive answer to Question 1.

Let  $\A$ and $\B$ be two abelian categories,  $F : \mathcal{A} \rightarrow \mathcal{B}$ a right exact additive functor, and $G:\B\to \A$ a left exact additive functor. Recall from \cite{FGR75} that there are two abelian categories $(F \downarrow \B)$ and $(\A \uparrow G)$ called comma categories. When $F$ is left adjoint to $G$, $(F \downarrow \B)$ is isomorphic to $(\A \uparrow G)$. Throughout this introduction, we assume $\A$ and $\B$ have enough projective and injective objects, direct sums and  products, and that they are exact (though these assumptions are not all needed in our main results and are made only for the sake of a well-presented Introduction; see Setups I and II). If $(F,G)$ is an adjoint pair, both comma categories are identified, denoted as $\T$, and called the comma category (see Subsection \ref{Subsec: comma category} for more details). 

In this paper, we introduce the notions of totally compatible and totally cocompatible functors and bimodules, and then provide an answer to Question 1 for comma categories (Theorem \ref{G-proj/inj}) and, as an immediate consequence, for triangular matrix rings (Corollary \ref{G-proj/inj mod}).

\bigskip

\noindent {\bf Theorem A.} Assume $(F,G)$ is an adjoint pair. Then 
\begin{enumerate}
	\item  $\GP(\T)=\Phi(\GP(\A),\GP(\B))$  if and only if $F$ is totally compatible.
	\item $\GI(\T)=\Psi(\GI(\A),\GI(\B))$  if and only if $G$ is totally cocompatible.
\end{enumerate}

\bigskip

Note that, besides the obvious gain of generality, there is another advantage to working with comma categories. While it is not true that the opposite of
a module category is a module category, it is true that the opposite of a comma category is a comma category (see Lemma \ref{the comma-categ-op}). This observation  makes it easy to dualize the results concerning Gorenstein projective  and, in a sense, cuts the work in half.

\bigskip

Recall that the global dimension of a ring $R$ is defined as the supremum of the projective (resp., injective) dimensions of all $R$-modules. Among other things, it measures ‘how far’ a ring is from possessing certain homological properties. For instance, a ring is semisimple if and only if its global dimension is zero, and a commutative noetherian local ring is regular if and only if its global dimension is finite; see, e.g.  \cite[Sections 8.1 and 8.4]{Rot79}.

A characterization of when the global dimension of a triangular matrix ring is finite/bounded has been given by Palmer and Roos \cite[Theorem 5]{PR73}, using classical derived functors. Recently, Minamoto and Yamaura also provided a modern expression of a formula for the global dimension of a triangular matrix algebra (see \cite[Corollary 6.2]{MY20}). The problem of finding similar formulas in Gorenstein homological algebra remains unresolved. This leaves us with another natural open question:

\bigskip

 \textbf{Question 2:} Is there an exact formula for computing the global Gorenstein dimension of the triangular matrix ring $T$ using the global Gorenstein dimension of $A$ and $B$?

 \bigskip
Recall that by \cite[Theorem A.6]{CR22}, the global Gorenstein dimension of an abelian category $\C$ with enough projective and injective objects is defined as the supremum of the set of Gorenstein projective (resp., injective) dimensions of all objects in $\C$:
$$\glGdim(\C):=\sup\{\Gpd_\C(M)  | M \in \C\}=\sup\{\Gid_\C(M)  | M \in \C\}$$
where $\Gpd_\C(M)$ (resp., $\Gid_\C(M)$) denotes the Gorenstein projective (resp., injective) dimension of an object $M$ in $\C$.

 Our second aim in the paper is to study the behavior of the global Gorenstein dimension of the comma category. Our approach consists of using Theorem A to provide an estimation of the global Gorenstein dimension of the comma category $\T$ using the global Gorenstein dimensions of the categories $\A$ and $\B$ and the global Gorenstein dimensions of the functors $F:\A\to \B$ and  $G:\B\to \A$ defined, respectively, as:
 $$\GPD(F):=\sup\{\Gpd_\B(F(A))   \,|\,  A\in \GP(\A)  \},$$ 
 $$\GID(G):=\sup\{\Gid_\A(G(B))   \,|\,  B\in \GI(\B)  \}.$$
see Theorem \ref{global-estimation}.
   
\bigskip

\noindent {\bf Theorem B.} Assume $(F,G)$ is an adjoint pair. 
 \begin{enumerate}
     \item 	If $F$ is totally compatible, then $$\max\{\glGdim(\A),\glGdim(\B)\} $$ 
      $$\leq \glGdim(\T)
\leq $$   $$ \max\{\glGdim(\A)+\GPD(F)+1,\glGdim(\B)\}$$
\item If $G$ is totally cocompatible, then
$$\max\{\glGdim(\A),\glGdim(\B)\}  $$  $$\leq \glGdim(\T)\leq $$  $$ \max\{\glGdim(\A), \glGdim(\B) +\GID(G)+1\}.$$

 \end{enumerate}

 Consequently, if $F$ is totally compatible  or $G$ is totally cocompatible, then the comma category $\T$  is of finite global Gorenstein dimension if and only if so are $\A$ and $\B$.

\bigskip

The approach used to prove Theorem B applies not only to the (global) dimensions defined via Gorenstein objects but also to (global) dimensions with respect to any class of the form $\X=\Phi(\X_1,\X_2)$ (resp., $\Y=\Psi(\Y_1,\Y_2)$); see Theorems \ref{glresdim of comma} and \ref{glcoresdim of comma}.

Motivated by Question 2, we investigate when the estimation in Theorem B becomes an exact formula. This gives several situations in which we have a partial answer to this open question (see Propositions \ref{Formula 1}, \ref{Formula 2} and  \ref{Formula 3}).

\bigskip

\noindent {\bf Theorem C.} Assume $(F,G)$ is an adjoint pair.

\begin{enumerate}
    \item If  $F=0$ or $G=0$, then 		 
          $$\glGdim(\mathcal{T})=\glGdim(\A\times \B)= \max\{\glGdim (\A),\glGdim(\B)\}.$$
     \item  If $\A$ is Frobenius and $0\neq F:\A \to \B$ is weakly compatible, then 
             $$\glGdim(\T)=
			\max\{\GPD(F)+1,\glGdim(\B)\}.$$
        \item  If $\B$ is Frobenius and $0\neq G:\B \to \A$ is weakly cocompatible, then 
          $$ \glGdim(\T)= \max\{\glGdim(\A),\GID(G)+1\}.$$
         \item If $(F, G)$ is a Frobenius pair with $0\neq F:\A \to \B$ faithful, then
       $$ \glGdim(\T)=\max\{\glGdim(\A)+1,\glGdim(\B)\}.$$

    \item  If  $(F, G)$ is a Frobenius pair with $0\neq G:\B \to \A$ faithful, then 
         $$ \glGdim(\T)=\max\{\glGdim(\A),\glGdim(\B)+1\}.$$
 
\end{enumerate}

We refer the reader to Section~\ref{Sec. GPD} for the  definitions of Frobenius categories and functors.

\bigskip

Another category of interest that we deal with in this paper is  the $n$-morphism category $\Mor_n(\C)$, where $n\geq 2$ is an integer and $\C$ is an abelian category. This category arises naturally in the representation theory of Artin algebras. For $n=2$, this is nothing but the well-known morphism category. For an arbitrary $n$, it was first introduced by Zhang in \cite{Zha11} in module categories and further defined and investigated for abelian categories by Wang and Liu \cite{WL22} from the viewpoint of triangulated categories and Gorenstein homological algebra (see Example \ref{Mor-cate} for more details on this category). As the $n$-morphism category can be seen as a  comma category (Lemma \ref{iso of n-Mor}), a direct application of Theorem A gives a complete description of Gorenstein projective/injective objects in $\Mor_n(\C)$ via Zhang's monomorphism and epimorphism categories  (Corollary  \ref{G-proj/inj in Mor_n}), while an application of Theorem B gives an exact formula to compute the global Gorenstein dimension of $\Mor_n(\C)$ via the global Gorenstein dimension of $\C$ (Proposition \ref{Formula 4: Mor_n}).

\bigskip

\section{Preliminaries}

Throughout, all rings are assumed to be associative with identity,  and  \( R\text{-Mod} \) will denote the category of all (left) modules over a ring $R$. We adopt the convention that an \( R \)-module always refers to a left \( R \)-module. By a subcategory of an abelian category \( \mathcal{C} \), we always mean a full subcategory closed under isomorphisms. Any class of objects in \( \mathcal{C} \) will be regarded as a (full) subcategory. Finally, we assume that all functors  are additive.

\subsection{Comma categories} \label{Subsec: comma category}

Throughout the paper, we will assume that $\A$ and $\B$ are abelian categories, $F : \mathcal{A} \rightarrow \mathcal{B}$ a right exact  functor, and $G:\B\to \A$ a left exact functor. It follows from \cite{FGR75} that there are two abelian categories called comma categories, denoted  $(F \downarrow \B)$  and $(\A \uparrow \G)$. 

The comma category $(F \downarrow \mathcal{B})$ is defined as follows:  The objects are triples $\begin{psmallmatrix}
    M_1\\ M_2
\end{psmallmatrix}_{\varphi^M}$, with $M_1 \in \mathcal{A}$, $M_2 \in \mathcal{B}$ and $\varphi^M : F(M_1) \rightarrow M_2$ is a morphism in $\mathcal{B}$. A morphism $\begin{psmallmatrix}
    f_1\\ f_2
\end{psmallmatrix} : \begin{psmallmatrix}
    M_1\\ M_2
\end{psmallmatrix}_{\varphi^M}\rightarrow \begin{psmallmatrix}
    N_1\\ N_2
\end{psmallmatrix}_{\varphi^N}$ is given by two morphisms $f_1 : M_1 \rightarrow N_1$ in $\mathcal{A}$ and $f_2 : M_2 \rightarrow N_2$ in $\mathcal{B}$ such that the following diagram is commutative
\[ \xymatrixcolsep{2pc} \xymatrixrowsep{1pc} 
\xymatrix{
	F(M_1) \ar[r]^{F(f_1)} \ar[d]_{\varphi^M} & F(N_1) \ar[d]^{\varphi^N} \\
	M_2 \ar[r]^{f_2} & N_2
}
\]
Dually, we can define the comma category $(\A\uparrow G)$. The objects are triples $\begin{psmallmatrix}
    M_1\\ M_2
\end{psmallmatrix}_{\psi^M}$, with $M_1 \in \mathcal{A}$, $M_2 \in \mathcal{B}$ and $\psi^M : M_1 \rightarrow G(M_2)$ is a morphism in $\A$, while morphisms are pairs $(f_1,f_2)$ making the obvious square commute. If there is no possible confusion, we sometimes omit the morphisms $\varphi^M$ and $\psi^M$.

In the following two lemmas, we give some relationships needed in this paper between the comma categories $(F \downarrow \mathcal{B})$ and  $(\A\uparrow G)$.

\begin{lem} \label{the comma-categ}  Assume $(F,G)$  is an adjoint pair.  There is an isomorphism of categories
$$(F \downarrow \mathcal{B}) \cong  (\A\uparrow G).$$

\end{lem}
\proof For each object $A\in \A$ and $B\in \B$,  there is a natural isomorphism $$\tau:=\tau_{A,B}:\: \Hom_\B(F A,B) \longrightarrow \Hom_\A(A,GB).$$ 
So, the isomorphism $ (F \downarrow \mathcal{B}) \to (\A\uparrow G)$ of categories is given on objects and morphisms, respectively, by $$\begin{psmallmatrix}
    M_1\\ M_2
\end{psmallmatrix}_{\varphi^M}\mapsto \begin{psmallmatrix}
    M_1\\ M_2
\end{psmallmatrix}_{\tau(\varphi^M)},\quad \begin{psmallmatrix}
    f_1\\ f_2
\end{psmallmatrix}\mapsto \begin{psmallmatrix}
    f_1\\ f_2
\end{psmallmatrix}$$ 
with inverse $(\A\uparrow G) \to (F \downarrow \mathcal{B})$ given on objects and morphisms, respectively, by $$\begin{psmallmatrix}
    M_1\\ M_2
\end{psmallmatrix}_{\psi^M}\mapsto \begin{psmallmatrix}
    M_1\\ M_2
\end{psmallmatrix}_{\tau^{-1}(\psi^M)} ,\quad \begin{psmallmatrix}
    f_1\\ f_2
\end{psmallmatrix}\mapsto \begin{psmallmatrix}
    f_1\\ f_2
\end{psmallmatrix}$$
\cqfd

\begin{defn}Under the isomorphism of Lemma \ref{the comma-categ},  both categories are identified, denoted as $\T$, and referred to as  \textbf{the comma category}.
\end{defn}

Throughout, we denote the opposite category of a  category $\C$ by $\C^{op}$ and the opposite functor of a functor $H:\C\to \mathcal{D}$ by $H^{op}:\C^{op}\to \mathcal{D}^{op}$.

\begin{lem} \label{the comma-categ-op} There are isomorphisms of categories $$(\A \uparrow G)^{op} \cong (G^{op}\downarrow \A^{op}) \text{ and }  (F\downarrow \B)^{op}  \cong (\B ^{op} \uparrow F^{op}) $$
     
\end{lem}
\proof The first isomorphism $(\A \uparrow G)^{op} \to (G^{op}\downarrow \A^{op})$ is given on objects and morphisms, respectively, by $$\begin{psmallmatrix}
    M_1\\ M_2
\end{psmallmatrix}_{\psi^M}\mapsto \begin{psmallmatrix}
    M_2\\ M_1
\end{psmallmatrix}_{(\psi^M)^{op}} ,\quad \begin{psmallmatrix}
    f_1\\ f_2
\end{psmallmatrix}^{op}\mapsto \begin{psmallmatrix}
    f_2^{op}\\ f_1^{op}
\end{psmallmatrix}$$
with inverse  $$(G^{op}\downarrow \A^{op})\to (\A \uparrow G)^{op}, \begin{psmallmatrix}
    N_1\\ N_2
\end{psmallmatrix}_{\varphi^N}\mapsto \begin{psmallmatrix}
    N_2\\ N_1
\end{psmallmatrix}_{(\varphi^N)^{op}},\quad \begin{psmallmatrix}
    g_1^{op}\\ g_2^{op}
\end{psmallmatrix}\mapsto \begin{psmallmatrix}
    g_2\\ g_1
\end{psmallmatrix}^{op}$$

The second isomorphism is dual. However, using the above isomorphism, we can alternatively obtain:
$$ (\B ^{op} \uparrow F^{op}) \cong  ((\B ^{op} \uparrow F ^{op})^{op})^{op} \cong ((F^{op})^{op}\downarrow (\B^{op})^{op})^{op} \cong (F\downarrow \B)^{op}.$$
\cqfd

Given two classes of objects $\X_1\subseteq \A$ and $\X_2\subseteq \B$, we set 
$$\Phi(\X_1,\X_2):=\{M=\begin{psmallmatrix}
    M_1\\ M_2
\end{psmallmatrix}_{\varphi^M}\in (F \downarrow \B) \;|\; M_1\in \X_1, \Coker\varphi^M\in \X_2, \text{ $\varphi^M$ is monic}  \},$$
$$\Psi(\X_1,\X_2):=\{M=\begin{psmallmatrix}
    M_1\\ M_2
\end{psmallmatrix}_{\psi^M}  \in (\A\uparrow G) \;|\; \Ker \psi^M\in \X_1, M_2\in \X_2, \text{ $\psi^M$ is epic}  \}.$$

We now present several examples of comma categories that are generally important and particularly relevant to us in the context of the present paper.

\begin{exmp}[\textbf{Triangular matrix rings}] Let $A$ and $B$ be two rings, $U$ a $(B,A)$-bimodule, and $T = \begin{psmallmatrix} A & 0 \\ U & B \end{psmallmatrix}$ the triangular matrix ring. If we set $F = U \otimes_A - : A\text{-Mod} \rightarrow B\text{-Mod}$ and $G:=\Hom_B(U,-): B\text{-Mod} \rightarrow A\text{-Mod}$, it follows by \cite{FGR75,Gre82} that we have equivalences of categories
  $$T\text{-Mod} \simeq (F \downarrow B\text{-Mod}) \simeq  (A\text{-Mod} \uparrow G).$$

Throughout this paper,  $T=\begin{psmallmatrix}
A & 0\\ 
U & B  
\end{psmallmatrix}$ will always denote a triangular matrix ring, and $T$-Mod will be identified with  both categories $(F \downarrow B\text{-Mod})$  and $(A\text{-Mod} \uparrow G).$
\end{exmp}

\begin{exmp} [\textbf{n-morphism categories}] \label{Mor-cate} Let $\C$ be an abelian category. Denote by $\Mor(\C)$ the category of morphisms in $\C$: the objects are morphisms in $\C$ and the morphisms are given by commutative squares in $\C$. If $\id_\C$ is the identity functor, then the morphism category $\Mor(\mathcal{C})$ of $\mathcal{C}$  coincides with the comma category $(\id_\C \downarrow \mathcal{C})$. 

For a given integer $n\geq 1$, one may consider a more general situation. Recall from \cite{Zha11,WL22} that there is an abelian category $\Mor_n(\C)$, called the $n$-morphism category,  defined as follows: An object of the $n$-morphism category $\Mor_n(\C)$ is an $(n-1)$-tuple of composable morphisms, written as
 $$(X_i,\varphi_i^X):=X_1\stackrel{\varphi_1^X}{\longrightarrow} X_2 \stackrel{\varphi_2^X}{\longrightarrow}  X_3 \stackrel{\varphi_3^X}{\longrightarrow}  \cdots \stackrel{\varphi_{n-1}^X}{\longrightarrow}  X_n \text{ or } (X_i,\varphi_i^X):=\begin{psmallmatrix}
     X_1 \\ X_2\\ \vdots \\ X_n
 \end{psmallmatrix}_{(\varphi^X_i)}$$
A morphism $(X_i, \varphi_i^X) \to (Y_i, \varphi^Y_i)$ is $f=(f_i)$, where $f_i:X_i \to  Y_i$ is a morphism for each $i$ such that every square in the following diagram commutes
$$\xymatrix{  X_1\ar[r]^{\varphi_1^X} \ar[d]_{f_1}  & X_2  \ar[d]_{f_2}  \ar[r] ^{\varphi_2^X}  & \cdots \ar[r] ^{\varphi_{n-1}^X}    & X_n  \ar[d]_{f_n}   \\
            Y_1 \ar[r]^{\varphi_1^Y}  & Y_2 \ar[r]^{\varphi_2^Y} & \cdots \ar[r] ^{\varphi_{n-1}^Y}  & Y_n }$$

\end{exmp}

 The $1$-morphism category $\Mor_{1}(\C)$ is nothing but the category $\C$, while the $2$-morphism category $\Mor_{2}(\C)$ coincides with the morphism category $\Mor(\C)$. More generally, we have:

\begin{lem}\label{iso of n-Mor} For any integer $n\geq 2$, there are isomorphisms of categories
$$\Mor_{n}(\C) \cong (H_n \downarrow \Mor_{n-1}(\C)) \cong  (\Mor_{n-1}(\C)    \uparrow H_n),$$
where $H_n:\C \to \Mor_{n-1}(\C)$ is a functor defined on objects by $X\mapsto \begin{psmallmatrix}
     X \\ X\\ \vdots \\ X
 \end{psmallmatrix}_{(\id_X)_i}$

\end{lem}
\proof The isomorphism $\Mor_{n}(\C) \cong (H_n \downarrow \Mor_{n-1}(\C))$ is given by the functor defined on objects by 
$$\Mor_{n}(\C) \to  (H_n \downarrow \Mor_{n-1}(\C)), \begin{psmallmatrix}
     X_1 \\ X_2\\ \vdots \\ X_n
 \end{psmallmatrix}_{(\varphi^X_i)} \mapsto \begin{psmallmatrix}
    X_1  \\ \hspace{0.5cm}  \begin{psmallmatrix}
     X_2 \\ \vdots \\ X_n
 \end{psmallmatrix}_{(\varphi^X_i)}
\end{psmallmatrix}_{\varphi^X}$$
where $\varphi^X:=(\varphi_1^X,\varphi_2^X\varphi_1^X, \cdots, \varphi_{n-1}^X \cdots \varphi_2^X\varphi_1^X) $, with inverse
$$(H_n \downarrow \Mor_{n-1}(\C)) \to  \Mor_{n}(\C), \begin{psmallmatrix}
    X_1  \\ \hspace{0.5cm}  \begin{psmallmatrix}
     X_2 \\ \vdots \\ X_n
 \end{psmallmatrix}_{(\varphi^X_i)}
\end{psmallmatrix}_{\varphi^X}\mapsto  \begin{psmallmatrix}
     X_1 \\ X_2\\ \vdots \\ X_n
 \end{psmallmatrix}_{(\varphi^X_i)}$$
where $\varphi_1^X:X_1\to X_2$ is the first component of the morphism
$\varphi^X$.

It is immediate to verify that the two assignments are mutually inverse on objects and morphisms.  The proof of the second isomorphism is dual. \cqfd

Let $\X$ be a class of objects in $\C$. Recall that the monomorphism category $\Phi_n(\X)$ (also denoted as $S_n(\X)$) is defined to be the full subcategory of $\Mor_n(\C)$ consisting of all the objects $(X_i, \varphi_i^X)$ such that $X_i\in \X$ for $1\leq i\leq n$, $\varphi_i^X:X_i\to X_{i+1}$ is a monomorphism and $\Coker\varphi_i^X\in \X$ for $1\leq i\leq n-1$. Dually, the epimorphism category $\Psi_n(\X)$ (also denoted as $\F_n(\X)$) is defined to be the full subcategory of $\Mor_n(\C)$ consisting of all the objects $(X_i, \psi_i^X)$ such that  $X_i\in \X$ for each $1\leq i\leq n$, $\psi_i^X:X_i\to X_{i+1}$ is an epimorphism and $\Ker\psi_i^X\in \X$ for $1\leq i\leq n-1$. Throughout, we identify  $\Phi_{1}(\X)$ (resp., $\Psi_{1}(\X)$) with $\X$.

\begin{lem}\label{iso of Phi_n} Assume $\X$ is closed under extensions. For any integer $n\geq 2$, there are isomorphisms of categories
$$\Phi_n(\X)\cong \Phi(\X,\Phi_{n-1}(\X)) \text{ and }   \Psi_n(\X)\cong \Psi(\Psi_{n-1}(\X),\X).$$ 
 
\end{lem}

\proof It is enough to show that the first isomorphism of categories from Lemma \ref{iso of n-Mor} restricts to an isomorphism of categories $$\Phi_n(\X)\cong \Phi(\X,\Phi_{n-1}(\X)).$$

For an object  $X= \begin{psmallmatrix}
     X_1 \\ X_2\\ \vdots \\ X_n
 \end{psmallmatrix}_{(\varphi^X_i)} \in \Phi_n(\X)$, we claim that $$\begin{psmallmatrix}
    X_1  \\ \hspace{0.5cm}  \begin{psmallmatrix}
     X_2 \\ \vdots \\ X_n
 \end{psmallmatrix}_{(\varphi^X_i)}
\end{psmallmatrix}_{\varphi^X}\in \Phi(\X,\Phi_{n-1}(\X)).$$
 That is, $X_1\in \X$ and $\varphi^X$ is a monomorphism with $\Coker\varphi^X\in \Phi_{n-1}(\X)$.

As $X\in \Phi_n(\X)$, we have by definition that  $X_i\in \X$ for $1\leq i\leq n$, $\varphi_i^X:X_i\to X_{i+1}$ is a monomorphism and $\Coker\varphi_i^X\in \X$ for $1\leq i\leq n-1$. For each $i=2, \cdots, n-1$, the morphism $\varphi^X$ induces a commutative diagram
$$\xymatrix{  0 \ar[r] & X_1 \ar@{=}[d] \ar[r]^{\varphi^X_{i-1}\cdots \varphi^X_1}  & X_i  \ar[d]^{\varphi^X_i} \ar[r] & \Coker (\varphi^X_{i-1}\cdots \varphi^X_1) \ar[r] \ar@{-->}[d] & 0 \\
 0 \ar[r]   & X_1  \ar[r] ^{\varphi^X_i\cdots \varphi^X_1}         &X_{i+1} \ar[r] & \Coker (\varphi^X_{i}\cdots \varphi^X_1) \ar[r] & 0
}$$
It follows from the Snake Lemma that there exists a short exact sequence
$$0 \to \Coker (\varphi^X_{i-1}\cdots \varphi^X_1) \to \Coker (\varphi^X_{i}\cdots \varphi^X_1) \to \Coker \varphi^X_i  \to 0$$
By induction on $i$ we obtain $\Coker (\varphi^X_{i}\cdots \varphi^X_1)\in \X$ and hence $\Coker\varphi^X\in \Phi_{n-1}(\X)$. Note also that since each $\varphi_i^X$ is a monomorphism, every composite
$\varphi_i^X\cdots\varphi_1^X$ is a monomorphism, and hence
$\varphi^X$ is a monomorphism in $\Mor_{n-1}(\C)$.

Conversely,  for an object $X=\begin{psmallmatrix}
    X_1  \\ \hspace{0.5cm}  \begin{psmallmatrix}
     X_2 \\ \vdots \\ X_n
 \end{psmallmatrix}_{(\varphi^X_i)}
\end{psmallmatrix}_{\varphi^X}\in \Phi(\X,\Phi_{n-1}(\X)),$ we claim that 
  $ \begin{psmallmatrix}
     X_1 \\ X_2\\ \vdots \\ X_n
 \end{psmallmatrix}_{(\varphi^X_i)} \in \Phi_n(\X).$ 
That is, $X_i\in \X$ for $1\leq i\leq n$, $\varphi_i^X:X_i\to X_{i+1}$ is a monomorphism and $\Coker\varphi_i^X\in \X$ for $1\leq i\leq n-1$.  As $X\in \Phi(\X,\Phi_{n-1}(\X))$, $X_1\in \X$ and $\varphi^X$ is a monomorphism with $\Coker\varphi^X\in \Phi_{n-1}(\X)$.  Now, we have the following  commutative diagram
$$\xymatrixcolsep{3pc} \xymatrixrowsep{1.5pc} \xymatrix{  0 \ar[r] & X_1 \ar@{=}[d] \ar[r]^{\varphi^X_{1}}  & X_2  \ar[d]^{\varphi^X_2} \ar[r] & \Coker \varphi^X_{1}\ar[r] \ar@{-->}[d]^{\overline{\varphi^X_2}} & 0 \\
 0 \ar[r]   & X_1   \ar@{=}[d] \ar[r]   ^{\varphi^X_2 \varphi^X_1}         &X_{3} \ar[r] \ar[d]^{\varphi^X_3} & \Coker (\varphi^X_2 \varphi^X_1) \ar[r] \ar@{-->}[d]^{\overline{\varphi^X_3}}  & 0
\\
 0 \ar[r]   & X_1  \ar[r]^{\varphi^X_3 \varphi^X_2 \varphi^X_1}          &X_{4} \ar[r]  & \Coker (\varphi^X_{3} \varphi^X_2\varphi^X_1) \ar[r] & 0 \\
& \vdots & \vdots & \vdots  & }$$
which shows that $\Coker\varphi^X=\begin{psmallmatrix}
     \Coker \varphi^X_{1} \\ \Coker (\varphi^X_2 \varphi^X_1)\\ \vdots \\ \Coker (\varphi^X_{n-1}\cdots\varphi^X_2 \varphi^X_1)
\end{psmallmatrix}_{(\overline{\varphi^X_i})_i}$. Therefore, we get that $\Coker (\varphi^X_{i}\cdots \varphi^X_1)\in \X$, each  $\overline{\varphi^X_i}$ is a monomorphism and each $\Coker \overline{\varphi^X_i}\in \X$. Since $X_1\in \X$ and each $\Coker (\varphi^X_{i}\cdots \varphi^X_1)\in \X$, the fact that $\X$ is closed under extensions shows that $X_{i+1}\in \X$ for all $1\leq i\leq n-1$. Moreover, the Snake Lemma applied to the above diagram shows that 
$\Coker\varphi^X_i\cong \Coker\overline{\varphi^X_i} \in \X$ and that $\varphi^X_i: X_i\to X_{i+1}$ is a monomorphism for all $2\leq i\leq n-1$. \cqfd

\begin{exmp}[\textbf{Trivial extension}] \label{Triv-ext} Let $\C$ be an abelian category and $H:\C \to \C$ a right exact  functor. Recall from \cite{FGR75} that the (right) trivial extension of $\C$ by $H$, denoted $\C \ltimes H$, is the category defined by: the objects are pairs
$(X,f)$ with $X\in \C$ and 
$f:H(X)\to X$ a morphism in $\C$ such that the composition  $f\circ H(f)=0$.  A morphism 
$\gamma: (X,f)\to (Y,g)$ is a morphism $\alpha: X \to Y$ such that the diagram
\[ \xymatrixcolsep{1.3pc} \xymatrixrowsep{1.3pc}
\xymatrix{
	H(X) \ar[r]^{H(\alpha)} \ar[d]_{f} & H(Y) \ar[d]^{g} \\
	X \ar[r]^{\alpha} & Y
}
\]
commutes

If we take $\C=\A \times \B$, the functor $F:\A\to \B$ induces a functor $\widetilde{F} : \A \times \B \to  \A \times \B$ defined by 
$$\widetilde{F}(A,B)=(0,F(A)) \hspace{1cm} \widetilde{F}(\alpha,\beta)=(0,F(\alpha))$$
It follows by \cite{FGR75} that the comma category $(F \downarrow \B)$ is isomorphic to the trivial extension $(\A\times \B)\ltimes \widetilde{F}$. In this case, an object $\begin{psmallmatrix}
    A\\ B
\end{psmallmatrix}_f$ in $(F \downarrow \B)$ induces a morphism $(0,f): \widetilde{F}(A,B) \to (A,B)$ in $\A\times \B$ and hence an object $((A,B), (0,f))$ in  $(\A\times \B)\ltimes \widetilde{F}$. On the other hand, as $\widetilde{F}^2=0$,  any object $((A,B), (f,g))$ in $(\A\times \B)\ltimes \widetilde{F}$ induces a morphism $g:F(A)\to B$ in $\B$ and so an object $\begin{psmallmatrix}
    A \\ B
\end{psmallmatrix}_g$  in $(F \downarrow \B)$.

Given a left exact functor $H:\C\to \C$, one can define dually the (left) comma category $H\rtimes \C$. In the same way, we see that in case $\C=\A\times \B$, the functor $G$ induces a functor $\widetilde{G}: \A\times \B \to \A\times \B$ such that  the comma category $(\A \uparrow G)$ is isomorphic to the trivial extension $\widetilde{G}\rtimes (\A\times \B)$.

\end{exmp}

\noindent {\bf Canonical adjoint functors.} The comma categories admit several canonical adjoint functors  
$$(\pp,\qq), \quad (\qq,\hh), \quad (\cc,\zz), \quad (\zz,\kk)$$
relating them to the product category $\A\times \B$. These functors will play a central role in the study of projective, injective, and Gorenstein objects.
	
\begin{enumerate}
	\item (\textbf{The tensor functor}) $\pp: \A \times \B\rightarrow (F \downarrow \mathcal{B})$ 
	$$\pp(M_1,M_2)=\begin{psmallmatrix}
	M_1 \\ 
	F(M_1)\oplus M_2   
	\end{psmallmatrix}_{(1,0)} \text{ and } \pp(f_1,f_2)=\begin{psmallmatrix}
	f_1 \\ 
	F(f_1)\oplus f_2   
	\end{psmallmatrix}$$ 
	for each object $(M_1,M_2)$ and  morphism $(f_1,f_2)$ in $\A \times \B$. 
	
	\item (\textbf{The hom functor})  $\textbf{h}:\A \times \B\rightarrow (\A \uparrow \G)$: $$\textbf{h}(M_1,M_2)=\begin{psmallmatrix}
	M_1 \oplus G(M_2)\\ 
	M_2 
	\end{psmallmatrix}_{(0,1)} \text{ and } \textbf{h}(f_1,f_2)=\begin{psmallmatrix}
	f_1 \oplus G(f_2)\\ 
	f_2 
	\end{psmallmatrix}$$
    for each object $(M_1,M_2)$ and  morphism $(f_1,f_2)$ in $\A \times \B$. 

    	\item (\textbf{The underlying functor}) $\textbf{q}:\T\rightarrow \A \times \B$: 
	 $$\textbf{q}(M)
	=(M_1,M_2)\, \text{ and }\, \textbf{q}(f)
	=(f_1,f_2)$$
 for each object $M=\begin{psmallmatrix}
	M_1 \\ 
	M_2   
	\end{psmallmatrix}$ and  morphism  $f=\begin{psmallmatrix}
	f_1 \\ 
	f_2   
	\end{psmallmatrix}$ in $\T$, where $\T=(F \downarrow \mathcal{B})$ or $\T=(\A \uparrow \G)$.

	\item (\textbf{The cokernel functor})  $\cc:(F \downarrow \mathcal{B})\rightarrow \A \times \B$: $$\cc(M)
	=(M_1,\Coker\varphi^M) \text{ and } \cc(f)
	=(f_1,\overline{f}_2)$$ for each object $M=\begin{psmallmatrix}
	M_1 \\ 
	M_2   
	\end{psmallmatrix}_{\varphi^M}$ and morphism  $f=\begin{psmallmatrix}
	f_1 \\ 
	f_2   
	\end{psmallmatrix} :M\to N$  in $(F \downarrow \mathcal{B})$, where $\overline{f}_2:\Coker\varphi^M\to\Coker\varphi^N$ is the induced morphism.

    \item (\textbf{The kernel functor})  $\kk:(\A \uparrow \G)\rightarrow \A \times \B$: $$\kk(M)
	=(\Ker\psi^M, M_2) \text{ and } \kk(f)
	=(\underline{f_1},f_2)$$ 
for each object $M=\begin{psmallmatrix}
	M_1 \\ 
	M_2   
	\end{psmallmatrix}_{\psi^M}$ and morphism  $f=\begin{psmallmatrix}
	f_1 \\ 
	f_2   
	\end{psmallmatrix}:M\to N$  in $(\A \uparrow \G)$, where $\underline{f_1}:\Ker\psi^M\to\Ker\psi^N$ is the induced map.

    \item (\textbf{The zero functor}) $\zz$\;:\;$\A \times \B\rightarrow \T$  $$\zz(M_1,M_2)=\begin{psmallmatrix}
	M_1 \\ 
	M_2   
	\end{psmallmatrix}_0 \text{ and } \zz(f_1,f_2)=\begin{psmallmatrix}
	f_1 \\ 
	f_2   
	\end{psmallmatrix}$$ for any object $(M_1,M_2)$
	and  morphism $(f_1,f_2)$ in  $\A \times \B$, where $\T=(F \downarrow \mathcal{B})$ or $\T=(\A \uparrow \G)$.
\end{enumerate}

\begin{rem}\item 
\begin{enumerate}
\item Note that we have used the functor $\qq$ to denote two different functors: $\textbf{q}:(F \downarrow \mathcal{B})\rightarrow \A \times \B$ and $\textbf{q}\;:\;(\A \uparrow \G)\rightarrow \A \times \B$. This will not cause any problems, as it will be clear from the context which functor we mean. The same remark applies to the functor $\zz$.

\item Since $\qq$ is exact and $\pp$ (resp., $\hh$) is left (resp., right) adjoint to $\qq$, $\pp$ (resp., $\hh$) preserves projective (resp., injective) objects.

\end{enumerate}

\end{rem}

Denote by $L_i F$ (resp., $R^iG$)  the left  (resp., right) derived functors of $F$ (resp., $G$). 

\begin{lem}\label{Ext-iso} Let $n\geq 1$ be an integer. Then, we have the following natural isomorphisms:
	\begin{enumerate}
		
	\item Let $N=\begin{psmallmatrix}
	N_1\\N_2
	\end{psmallmatrix}_{\varphi^N}\in (F \downarrow \mathcal{B})$ and  $(M_1, M_2)\in\A\times  \B$. 	\begin{enumerate}
		    \item If ${\rm L}_{i}F(M_1)=0$ for all $1\leq i\leq n$, then $$\Ext^n_{(F \downarrow \mathcal{B})}(\begin{psmallmatrix}
		M_1\\ F(M_1)
		\end{psmallmatrix},N)\cong \Ext^n_\A(M_1,N_1).$$
		\item  $\Ext^n_{(F \downarrow \mathcal{B})}(\begin{psmallmatrix}
		0\\M_2
		\end{psmallmatrix},N) \cong \Ext^n_\B(M_2,N_2).$
		\end{enumerate}
        
	\item Let $M=\begin{psmallmatrix}
	M_1\\M_2
	\end{psmallmatrix}_{\psi^M}\in (\A \uparrow \G)$ and  $(N_1,N_2)\in  \A\times  \B$. 
    \begin{enumerate}
		    \item $\Ext^n_{(\A \uparrow \G)}(M,\begin{psmallmatrix}
		N_1\\0
		\end{psmallmatrix})\cong \Ext^n_\A(M_1,N_1).$
		\item If ${\rm R}^{i}G(N_2)=0$ for all $1\leq i\leq n$, then $$\Ext^n_{(\A \uparrow \G)}(M,\begin{psmallmatrix}
		G(N_2)\\N_2
		\end{psmallmatrix})\cong \Ext^n_\B(M_2,N_2).$$
		\end{enumerate}
	\end{enumerate}
\end{lem}

\proof  The proof is similar to that of triangular matrix rings; see \cite[Lemma 3.2]{Mao20} or  \cite[Lemma 3]{BEGO21}.  \cqfd

The following lemma is fundamental in this paper. Over triangular matrix rings, it follows from \cite[Theorem 3.1]{HV00} and \cite[Proposition 5.1]{HV99}.

\begin{lem} The following two equalities hold:\label{proj-inj in comma} \item 
	\begin{enumerate}
		\item $\P(F \downarrow \mathcal{B})=\Phi(\P(\A),\P(\B))$

		\item  $\I(\A \uparrow G)=\Psi(\I(\A),\I(\B))$   	\end{enumerate}
\end{lem}
\proof  Keeping in mind Example \ref{Triv-ext}, it follows from \cite[Corollary 1.6(c)]{FGR75} that an object $P=\begin{psmallmatrix}
    P_1 \\ P_2
\end{psmallmatrix}_{\varphi^P}$ in $(F \downarrow \mathcal{B})$ is projective if and only if $\cc(P)=(P_1,\Coker\varphi^P)$ is projective in $\A\times \B$ and $P$ is isomorphic to $\pp \cc(P)=\begin{psmallmatrix}
    P_1 \\ F(P_1)\oplus \Coker\varphi^P
\end{psmallmatrix}$. This is equivalent to saying that $P$  is projective if and only if $P_1$ is projective in $\A$, $\Coker\varphi^P$ is projective in $\B$ and $\varphi^P:F(P_1)\to P_2$ is a (split) monomorphism. By \cite[Corollary 1.6(d)]{FGR75}, a dual description of injective objects in $(\A \uparrow G)$ holds. \cqfd

\begin{lem}  The following assertions hold:\label{enough-proj/inj}\item
\begin{enumerate}
    \item If $\A$ and $\B$ have enough projective objects, so does $(F \downarrow \mathcal{B})$.
    \item If $\A$ and $\B$ have enough injective objects, so does $(\A \uparrow \G)$.
    \item Assume $(F,G)$ is an adjoint pair. If $\A$ and $\B$ have enough projective and injective objects, so does $\T$.
\end{enumerate}

\end{lem}
\proof   Let $M$ be an object in $(F \downarrow \mathcal{B})$. Since $\A$ and $\B$ have enough projective objects, there exist two epimorphisms   $f_1:P_1 \to M_1$ and $g_2:P_2 \to M_2$
in $\A$ and $\B$, respectively, where $P_1$ and $P_2$ are projectives. Since $F$ is right exact, we can construct a commutative diagram in $\B$
$$\xymatrixcolsep{1.3pc} \xymatrixrowsep{1.3pc} \xymatrix{   F(P_1) \ar[r]^{F(f_1)} \ar[d]^{\lambda_1}& F(M_1) \ar[r] \ar[d]^{\varphi^M} & 0   \\
F(P_1) \oplus P_2 \ar@{-->}[r]^{f_2} \ar[d]^{\pi_2}  & M_2 \ar[d] &    \\
 P_2 \ar@{-->}[r] \ar@/^1pc/[u]^{\lambda_2} \ar[ur]^{g_2}  & \Coker\varphi^M &   }$$
where $\lambda_i$ and $\pi_i$ are  the canonical injections and projections, respectively, and with $$f_2:F(P_1) \oplus P_2\to M_2$$ the unique morphism induced by $\varphi^M F(f_1)$ and $g_2$. Since $g_2=f_2 \circ\lambda_2$ is an epimorphism, so is $f_2$. Therefore, as $\pp$ preserves projective objects and both $f_1$ and $f_2$ are epimorphisms, we obtain an epimorphism $\pp(P_1,P_2) \to M$ with $\pp(P_1,P_2)$ projective. 

The statement for injective objects follows dually. \cqfd

 \begin{cor}\label{proj-inj in n-Mor}  For any integer $n\geq 1$, we have the following assertions: 
 \begin{enumerate}
		\item $\P(\Mor_n(\C))=\Phi_n(\P(\C)).$
		
		\item  $\I(\Mor_n(\C))=\Psi_n(\I(\C)).$   	\end{enumerate}
     
 \end{cor}
\proof By Lemmas \ref{iso of n-Mor} and \ref{iso of Phi_n},  we have isomorphisms of categories
$$\Mor_{n}(\C) \cong (H_n \downarrow \Mor_{n-1}(\C)) \text{ and }\Phi_n(\P(\C))\cong \Phi(\P(\C),\Phi_{n-1}(\P(\C))) \quad\forall n\geq 2,$$
with $\Mor_{1}(\C)=\C$ and $\Phi_{1}(\P(\C))=\P(\C)$. Hence, we obtain our desired description of projective objects in $\Mor_n(\C)$, using Lemma \ref{proj-inj in comma} and induction. 

The statement for injective objects follows dually. \cqfd

\subsection{$\X$-(co)resolution dimensions}

 Given a class $\X$ of objects of $\C$, recall that an $\X$-resolution of an object $A\in \C$ is an exact sequence 
$\cdots \to X_1\to X_0\to A\to 0$
with each $X_i\in\X$. An object $A\in\C$ is said to have $\X$-resolution dimension at most an integer $n \geq  0$, $\resdim_\X(A)\leq  n$, if $A$ has a finite $\X$-resolution: 
$$0\to X_n\to \cdots \to X_1\to X_0\to A\to 0.$$
If $n$ is the least non-negative integer for which such a sequence exists, then its $\X$-resolution dimension is precisely $n$, and if there is no such $n$ then we define its $\X$-resolution dimension as $\infty$.

Given a class $\mathcal{Z}$ of objects in $\C$, the $\X$-resolution dimension of $\mathcal{Z}$ is defined as: 
	$$\resdim_\X(\mathcal{Z})=\sup\{\resdim_\X(Z)\mid Z\in\mathcal{Z}\}.$$
In particular, we define the global $\X$-resolution dimension of $\C$ as the $\X$-resolution dimension of $\C$, that is, $$\glresdim_\X(\C):=\resdim_\X(\C)=\sup\{\resdim_\X(Z)\mid Z\in \C\}.$$

One can define $\X$-coresolutions and (global) $\X$-coresolution dimension dually.

Recall that a class of objects in $\C$ is projectively resolving if it is closed under extensions and  kernels of epimorphisms and contains all projective objects of $\C$. Injectively coresolving classes are defined dually.

\begin{lem}[\cite{Hua22}, Theorem 3.2] \label{resdim in ses} Let  $0\to N\to M\to L\to 0$ be a short exact sequence in an abelian category $\C$. Assume $\C$ has enough projective objects, and $\X$ a projectively resolving class which is closed under direct summands.
        
        \begin{enumerate}
			\item $\resdim_\X(M)\leq\sup\{\resdim_\X(N),\resdim_\X(L)\}.$\\ The equality holds if $\resdim_\X(N)+1\neq \resdim_\X(L)$. 
			\item $\resdim_\X(N)\leq\sup\{\resdim_\X(M),\resdim_\X(L)-1\}.$
            \\ The equality holds if $\resdim_\X(M)\neq \resdim_\X(L)$. 
			\item $\resdim_\X(L)\leq \sup\{\resdim_\X(N)+1,\resdim_\X(M)\}.$
            \\ The equality holds if $\resdim_\X(N)\neq \resdim_\X(M)$. 
		\end{enumerate}
	
\end{lem}

Let $\C$ be an abelian category with direct sums and they are exact. It is standard, by Eilenberg's swindle (see, e.g., \cite[Proposition 1.4]{Hol04} or \cite[Proposition 9.2.2]{Per16}), that a class of objects in $\C$ that is closed under direct sums and kernels of epimorphisms is also closed under direct summands. The following lemma follows from \cite[Lemma~5.1(1)]{Elm26}, in which the hypothesis that $\mathcal{X}$ is
projectively resolving was not assumed.

\begin{lem} [\cite{Elm26}, Lemma 5.1(1)]\label{X_n-closed-dir-sums/summ} Let $\C$ be an abelian category with enough projective objects and direct sums and they are exact, and let $\X$ be a subcategory of $\C$ that is projectively resolving and closed under direct sums. For any family $(M_i)_{i \in I}$ of objects  in $\C$, we have
		$$
		\resdim_\X\left(\bigoplus_{i \in I} M_i\right) = \sup\{\resdim_\X(M_i) \mid i \in I\}.
		$$
\end{lem}

Note that Lemmas \ref{resdim in ses}  and \ref{X_n-closed-dir-sums/summ}  have their dual counterparts (see \cite[Theorem 3.9]{Hua22} and  \cite[Lemma 5.1(2)]{Elm26}).

\subsection{Gorenstein projective and injective objects and dimension}

A complex $\PP$ of projective
objects in an abelian category $\C$ is called totally acyclic (or complete projective resolution) if it is acyclic and the induced complex $\Hom_\C(\PP,Q)$ is acyclic for every projective $Q\in\C$.  An object \(M \in \C\) is called Gorenstein projective if there exists a totally acyclic complex of projectives $\PP$ such that \(M \cong \Ker(P_{0}\to P_{-1})\). We denote the class of Gorenstein projective objects in $\C$ as $\GP(\C)$.

Dually, a complex $\EE$ of injective objects in $\C$ is called totally acyclic (or complete injective resolution) if it is acyclic and the induced complex \(\Hom_{\mathcal{C}}(E, \EE)\) is acyclic for every injective object \(E \in \C\).  An object \(N \in \C\) is called Gorenstein injective if there exists a totally acyclic complex of injectives $\EE$ such that \(N \cong \operatorname{Ker}(E_0\to E_{-1})\). We denote the class of Gorenstein injective objects in $\C$ as $\GI(\C)$.

Let $M$ be an object in $\C$. The Gorenstein projective dimension of $M$, $\Gpd_\C(M)$,  is defined as the $\GP(\C)$-resolution dimension of $M$. Dually, the Gorenstein injective dimension of $M$, $\Gid_\C(M)$,  is defined as the $\GI(\C)$-coresolution dimension of $M$. 

Recall also that the global Gorenstein projective dimension of $\C$, $\GPD(\C)$,  is defined as the global $\GP(\C)$-resolution dimension of $\C$, that is, 
$$\GPD(\C):=\glresdim_{\GP(\C)}(\C)= \sup\{\Gpd_\C(M)|\,M\in\C\}.$$
Dually, the global Gorenstein injective dimension of $\C$ is defined as 
$$\GID(\C):=\glcoresdim_{\GI(\C)}(\C)= \sup\{\Gid_\C(M)|\,M\in\C\}.$$
A basic result in Gorenstein homological algebra states that for an abelian category $\C$ with enough projective and injective objects, its global Gorenstein projective dimension coincides with its global Gorenstein injective dimension: 
$\GPD(\C)=\GID(\C).$
 This common value is called the global Gorenstein dimension of $\C$, denoted by $\glGdim(\C)$. This fact was first proved by Enochs and Jenda \cite[Section 12.3]{EJ00} for Gorenstein rings. It was later extended in \cite[Theorem 1.1]{BM10} to arbitrary rings  and recently in \cite{CR22} to abelian categories.

The following follows from \cite[Theorem A.6]{CR22}.
\begin{lem} \label{glGdim-charac} Let $\C$ be an abelian category  with enough projectives and injectives. Then
    $$\glGdim(\C)=\sup\{ \pd_\C(\I(\C)),\id_\C(\P(\C))\}.$$
\end{lem}

\begin{lem}\label{tot-acyc=acyc} Let $\C$ be an abelian category.
    \begin{enumerate}
        \item  If $\id_\C(\P(\C))<\infty$, then any acyclic complex of projective objects in $\C$ is totally acyclic.
    \item If $\pd_\C(\I(\C))<\infty$, then any acyclic complex of injective objects in $\C$ is totally acyclic.
        
    \end{enumerate}
\end{lem}
\proof This follows by  a standard argument using induction.   \cqfd

 \section{Gorenstein projective and injective objects in comma categories} \label{Sec. GP/GI}

 In this section, we prove Theorem A from the Introduction and give some applications.

\bigskip

\textbf{Setup I.} For the rest of the paper, we assume $\A$ and $\B$ have enough projective (resp., injective) objects whenever we are dealing with the comma category $(F\downarrow \B)$ (resp., $(\A \uparrow G)$). It follows from Lemma  \ref{enough-proj/inj} that $(F\downarrow \B)$ (resp., $(\A \uparrow G)$) has enough projective (resp., injective) objects.

\bigskip

In \cite[Definition 1.1]{Zha13}, Zhang introduced compatible bimodules in order to characterize Gorenstein projective modules over a triangular matrix algebra $T=\begin{psmallmatrix}
    A & 0 \\ U & B
\end{psmallmatrix}$. This was recently extended to comma categories independently by Hu and Zhu \cite{HZ22}, and Peng, Zhu, and Huang \cite{PZH22}. The first assumption in Zhang's definition of compatible bimodules is that all acyclic complexes of projective $A$-modules must be $(U\otimes_A-)$-exact. However, this condition seems to be stronger than necessary.  Our idea here is to replace acyclic complexes of projective objects with totally acyclic complexes. This leads to weaker conditions which turn out to be sufficient for the characterization of Gorenstein  objects in the comma categories.

\begin{defn} \label{P-tot-comp} A right exact functor $F: \A\to \B$ is \textbf{totally compatible} if it satisfies the following conditions:
	
	\begin{enumerate}
		
		\item[(P1)]  The complex $F(\PP)$ is acyclic for any totally acyclic complex $\textbf{P}$  of projective objects in $\A$.
		\item[(P2)]  The complex $\Hom_\B(\QQ,F(P))$ is acyclic for any totally acyclic complex $\QQ$ of projectives in $\B$ and any projective object $P$ in  $\A$.
		\item[(P3)]  The complex $F(\PP)$ is acyclic for any totally acyclic complex  $\begin{psmallmatrix}
\textbf{P}\\\textbf{Q}
		\end{psmallmatrix}$ of projective objects in $(F\downarrow \B)$.
	\end{enumerate}

\end{defn}

 Recall from \cite[Definition 4.1]{HZ22} that a right exact functor $F: \A\to \B$ is called weakly compatible if it satisfies conditions (P1) and (P2), while it is called compatible (or perfect in the sense of \cite{PZH22}) if it satisfies (P2) together with the following condition:
  \begin{enumerate}
	      \item[(P4)] $F(\PP)$ is exact for any exact complex of projective objects in $\A$.
	  \end{enumerate}
	
\begin{rem} We note that the results and definitions given in this paper always have their dual counterparts, as they can be obtained simply by passing to the opposite category (see Lemma \ref{the comma-categ-op}). Thus, in this paper, we only prove the statements concerning the comma category $(F\downarrow \B)$.
\end{rem}
 
\begin{defn}\label{I-tot-comp} A left exact functor  $G:\B \to \A$ is \textbf{totally cocompatible} if it satisfies the following conditions:
	\begin{enumerate}
		\item[(E1)] $G(\textbf{E}_2)$ is exact for any totally acyclic complex $\textbf{E}_2$  of injective objects in $\B$.
		\item[(E2)] $\Hom_\A(G(E),\textbf{E}_1)$ is exact for any totally acyclic complex $\textbf{E}_1$ of injective objects in $\A$ and any injective object $E$ in $\B$.
		\item[(E3)] $G(\textbf{E}_2)$ is exact for any totally acyclic complex $\begin{psmallmatrix}
		\EE_1\\ \EE_2
		\end{psmallmatrix}$ of injective objects in $(\A\uparrow G)$.
	\end{enumerate}	
\end{defn}

Following \cite[Definition 5.1]{YHW24}, the functor $G$ is called weakly cocompatible if it  satisfies $(E1)$ and $(E2)$, while it is called cocompatible if it satisfies $(E2)$ and the following assertion:
  \begin{enumerate}
	      \item[(E4)] $G(\textbf{E}_2)$ is exact for any acyclic complex $\textbf{E}_2$ of injective objects in $\B$.
	  \end{enumerate}

\begin{rem}\label{rem:tot-comp=comp} By \cite[Example 4.3]{HZ22}, a weakly compatible functor need not be compatible.
In general, we have the inclusions:
\begin{enumerate}
    \item[(i)]  $\text{\{compatible\}} \subseteq \text{\{totally compatible\}} \subseteq \text{\{weakly compatible\}}, $ \\
 with equalities when every projective $A\in \A$ has finite injective dimension.
     \item[(ii)] $\text{\{cocompatible\}} \subseteq \text{\{totally cocompatible\}} \subseteq \text{\{weakly cocompatible\}}$, \\
 with equalities when every injective $B\in \B$ has finite projective dimension.
\end{enumerate}
For example, the equalities in (i) (resp., (ii)) hold whenever $\A$ (resp., $\B$) has enough projective and injective objects and finite global Gorenstein dimension (see Lemmas \ref{glGdim-charac} and \ref{tot-acyc=acyc}).

Note that conditions (P1) and (P3) appear to be closely related. However, we do not know of any functor that satisfies (P1) but not (P3), or vice versa.
\end{rem}

Unlike projective (resp., injective) objects, Gorenstein projective (resp., Gorenstein injective) objects are not preserved by the functor $\pp$ (resp., $\h$). The following lemma characterizes precisely when this occurs.

\begin{lem}\label{Goren-preserv}  The following assertions hold.
	\begin{enumerate}
	    \item $\pp$ preserves Gorenstein projectives if and only if $F$ is weakly compatible. 
        \item $\hh$ preserves Gorenstein injectives if and only if $G$ is weakly cocompatible. 
	\end{enumerate}
\end{lem}
\proof The first claim follows from \cite[Lemma 4.6]{HZ22}, as $\pp\cong \pp(0,-)\oplus \pp(-,0)$ and the class of Gorenstein projective objects is closed under finite direct sums. The second assertion is dual.
\cqfd

We are now ready to prove Theorem A from the Introduction.
\begin{thm} \label{G-proj/inj} The following assertions hold.
	\begin{enumerate}
	\item  $\GP(F \downarrow \B)=\Phi(\GP(\A),\GP(\B))$  if and only if $F$ is totally compatible.
	\item $\GI(\A\uparrow G)=\Psi(\GI(\A),\GI(\B))$  if and only if $G$ is totally cocompatible.
	\end{enumerate}
\end{thm}
	
\proof  Throughout this proof, we set $\T=(F \downarrow \B)$.

$(\Leftarrow )$ Let $M=\begin{psmallmatrix}
M_1\\M_2
\end{psmallmatrix}_{\varphi^{M}}$ be an object in $\T$.

$(\subseteq)$ Assume that $M\in \GP(\T)$. By definition, there exists a totally acyclic complex of projective objects  in $\T$
$$\PP=\cdots\to P_1\to P_0\to P_{-1}\to P_{-2}\to \cdots$$
where $P_i=\begin{psmallmatrix}
P_{i,1}\\P_{i,2}
\end{psmallmatrix}_{\varphi^{P_i}}$ and such that $M \cong \Ker(P_0\to P_{-1}).$ Then, we get an acyclic complex in $\A$ 
$$\textbf{P}_1=\cdots\to P_{11}\to  P_{01}\to  P_{-11}\to P_{-21}\to \cdots$$
with each $P_{i1}$ projective by Lemma \ref{proj-inj in comma}, and such that  $M_1\cong \Ker(P_{01}\to 
 P_{-11}).$ 
 
 $\bullet$ $\varphi^M$ is a monomorphism: By  (P3), the complex $F(\PP_1)$ is exact. In particular, the short exact sequence
    $$F(\mathcal{E}_1): 0\to F(M_1)\to F(P_{01}) \to F(K_{-11})\to 0$$
  is exact, where $K_{-11}:=Z_{-1}(\PP_1)$. If $\iota_1:M_1\rightarrow P_{01}$ and $\iota_2:M_2\rightarrow P_{02}$ are the canonical monomorphisms, then $F(\iota_1)$ is a monomorphism, and from the following commutative diagram:
$$ \xymatrixcolsep{1.3pc} \xymatrixrowsep{1.3pc} \xymatrix{
 0\ar[r]	& F(M_1)\ar[r]^{F(\iota_1)} \ar[d]_{\varphi^{M}}  & F(P_{01})\ar[d]_{\varphi^{P_{0}}}\\ 
	& M_2  \ar[r]^{\iota_2}                 & P_{02},  
}  $$ 
we get that $\varphi^M$ is also a monomorphism.

$\bullet$ $\Coker\varphi^M$ is Gorenstein projective: Using again the fact that the complex $F(\textbf{P}_1)$ is exact and each $\varphi^{P^i}$ is a monomorphism, we obtain a commutative diagram with exact columns 
{\small $$\xymatrixcolsep{1.3pc} \xymatrixrowsep{1.3pc}\xymatrix{0 \ar[d]  & 0\ar[d]  & 0\ar[d]  &0\ar[d]  & 0 \ar[d]  \\
F(\textbf{P}_1):\ar[d]_{\varphi^{\textbf{P}}}  \cdots \ar[r]& F(P_{11})\ar[r] \ar[d]_{\varphi^{P^1}} & F(P_{01})\ar[d]_{\varphi^{P^0}} \ar[r]  &  F(P_{-11}) \ar[d]_{\varphi^{P^{-1}}} \ar[r]&  F(P_{-21})\ar[d]_{\varphi^{P^{-2}}} \ar[r] &  \cdots   \\
	  \textbf{P}_2: \ar[d]  \cdots \ar[r]& P_{12}\ar[r] \ar[d]& P_{02}\ar[d] \ar[r] & P_{-12}\ar[d] \ar[r] &  P_{-22} \ar[d]  \ar[r] & \cdots\\ 
	  \Coker \varphi^{\textbf{P}}: \ar[d]  \cdots\ar@{-->}[r] & \Coker \varphi^{P_1} \ar@{-->}[r] \ar[d]& \Coker \varphi^{P_0} \ar@{-->}[r] \ar[d] &  \Coker \varphi^{P_{-1}}  \ar[d] \ar[r]& \Coker \varphi^{P_{-2}} \ar@{-->}[r] \ar[d] & \cdots\\
	   0  &0 &0  &0 &0  & }$$} Since the first and second complexes are exact, we get an acyclic complex in $\B$
$$\Coker \varphi^{\textbf{P}}:\;\cdots\to \Coker \varphi^{P_1} \to \Coker \varphi^{P_0}\to  \Coker \varphi^{P_{-1}}  \to \Coker \varphi^{P_{-2}}\to  \cdots,$$
where each $\Coker \varphi^{P_{i}} $ is projective by Lemma \ref{proj-inj in comma}. Now, let $Q_2$ be a projective object in $\B$. Since $(\cc,\zz)$ is an adjoint pair, we obtain for each $i$ a natural isomorphism
\begin{eqnarray*}
    \Hom_ \B(\Coker\varphi^{P_i},Q_2) &\cong& \Hom_{\A \times \B}(\cc(P^i),(0,Q_2)  \\ &\cong &\Hom_{\T}( P_i,\zz(0,Q_2))=\Hom_{\T}(P_i, \begin{psmallmatrix}
    0 \\Q_2
\end{psmallmatrix}),
\end{eqnarray*}
which induces an isomorphism of complexes
$$ \Hom_ \B(\Coker\varphi^{\PP},Q_2)\cong \Hom_{\T}(\PP, \begin{psmallmatrix}
    0 \\Q_2
\end{psmallmatrix}).$$
By Lemma \ref{proj-inj in comma}, the object  $\begin{psmallmatrix}
    0 \\Q_2
\end{psmallmatrix}$ is projective. Hence the complex $\Hom_{\T}( \PP, \begin{psmallmatrix}
    0 \\Q_2
\end{psmallmatrix})$ is exact and so is the complex $\Hom_ \B(\Coker\varphi^{P},Q_2) $. Thus $\Coker \varphi^M$ is Gorenstein projective.

$\bullet$ $M_1$ is Gorenstein projective: We only need to show that $\Hom_\A(\PP_1,Q_1)$ is an acyclic complex for any  projective object $Q_1$ in $\A$. 

Applying the functor $\Hom_\T(\PP,-)$ to the following short exact sequence in $\T$
$$0\to \begin{psmallmatrix}
     0 \\ F(Q_1)
\end{psmallmatrix}\to \begin{psmallmatrix}
   Q _1 \\ F(Q_1)
\end{psmallmatrix}\to \begin{psmallmatrix}
    Q_1 \\ 0
\end{psmallmatrix}\to 0,$$
we get a short exact sequence of complexes
$$0\to \Hom_\T( \PP,\begin{psmallmatrix}
     0 \\ F(Q_1)
\end{psmallmatrix})\to \Hom_\T( \PP,\begin{psmallmatrix}
     Q_1 \\ F(Q_1)
\end{psmallmatrix})\to \Hom_\T( \PP,\begin{psmallmatrix}
     Q_1 \\ 0
\end{psmallmatrix})\to 0.$$
For each $i$, there exists a natural isomorphism of abelian groups
$$\Hom_\A(P_{i1},
     Q_1)\cong \Hom_{\T}( P_i,\begin{psmallmatrix}
     Q_1 \\ 0
\end{psmallmatrix}),$$ 
which  induces an isomorphism of complexes
   $$\Hom_\A(\PP_1,
     Q_1)\cong \Hom_{\T}( \PP,\begin{psmallmatrix}
     Q_1 \\ 0
\end{psmallmatrix}).$$ 
So the complex $\Hom_\A( \PP_1,
     Q_1)$ is acyclic if the complexes $$\Hom_\T( \PP,\begin{psmallmatrix}
     0 \\ F(Q_1)
\end{psmallmatrix}) \text{ and }\Hom_\T( \PP,\begin{psmallmatrix}
     Q_1 \\ F(Q_1)
\end{psmallmatrix})$$ are exact as well. 

The complex $\Hom_\T( \PP,\begin{psmallmatrix}
     Q_1 \\ F(Q_1)
\end{psmallmatrix})$ is exact, as $\begin{psmallmatrix}
     Q_1 \\ F(Q_1)
\end{psmallmatrix}$ is projective by Lemma \ref{proj-inj in comma}.

Using again the adjointness $(\cc,\zz)$, we obtain natural isomorphisms of complexes
\begin{eqnarray*}
    \Hom_{\T}( \PP,\begin{psmallmatrix}
     0 \\ F(Q_1)
\end{psmallmatrix})&=& \Hom_{\T}( \PP,\zz( 0 , F(Q_1))  \\ &\cong&   \Hom_{\A\times\B}(\cc(\PP),( 0 , F(Q_1)))  \\
&\cong& \Hom_{\B}(\Coker\varphi^{\PP} ,F(Q_1)).
\end{eqnarray*}
Thus $\Hom_{\T}( \PP,\begin{psmallmatrix}
     0 \\ F(Q_1)
\end{psmallmatrix})$ is exact by (P2), as $\Coker\varphi^{\PP}$ is totally acyclic.  Therefore, the complex $\Hom_\A(\PP_1,
     Q_1)$ is acyclic, and hence $M_1$ is Gorenstein projective.

\bigskip

$(\supseteq)$ Assume that $ M\in \Phi(\GP(\A),\GP(\B))$. By Lemma \ref{Goren-preserv}, the objects $\pp(M_1,0)$ and $\pp(0,\Coker\varphi^M)$ are Gorenstein projective. As the class of Gorenstein projective  objects is closed under extensions, it follows from the short exact sequence 
   $$0\to \pp(M_1,0) \to M\to \pp(0,\Coker\varphi^M)\to 0$$
that $M$ is Gorenstein projective.

$(\Rightarrow)$ The inclusion $\Phi(\GP(\A),\GP(\B)) \subseteq \GP(\T)$ shows that $\pp$ preserves Gorenstein projective objects, so, $F$ is weakly compatible by Lemma \ref{Goren-preserv}.

 It remains to show ${\rm (P3)}$ of Definition \ref{P-tot-comp}, that is, if $\PP=\begin{psmallmatrix}
\PP_1\\ \PP_2
\end{psmallmatrix}$ is a totally acyclic complex of projective objects in $\T$, then $F(\PP_1)$ is exact. Equivalently, we may show that for each short exact sequence
  $$0\to K_{k1}\stackrel{\iota_{k1}}{\to } P_{k1}\to K_{k-1,1}\to 0,$$
where $K_k=Z_k(\PP)=\begin{psmallmatrix}
K_{k1}\\K_{k2}
\end{psmallmatrix}_{\varphi^{K_k}}$, the short exact sequence 
$$0\to F(K_{k1})\stackrel{F(\iota_{k1})}{\to } F(P_{k1})\to F(K_{k-1,1})\to 0$$
is exact. Since $F$ is right exact, it suffices to show that $F(\iota_{k1})$ is a monomorphism.

By definition, each $K_k$ is Gorenstein projective. Therefore, each $\varphi^{K_k}$ is a monomorphism. It follows from the commutativity of the following diagram
$$  \xymatrixcolsep{1.3pc} \xymatrixrowsep{1.3pc} \xymatrix{
		&F(K_{k1})\ar[r]^{ F(\iota_{k1})} \ar[d]_{\varphi^{K_k}} \ar[r]  &F(P_{k1})\ar[d]_{\varphi^{P_k}}  \ar[r]  &F(K_{k-1,1})\ar[d]_{\varphi^{K_{k-1}}} \ar[r] & 0 \\   	0\ar[r]
	& K_{k2}  \ar[r]^{\iota_2}                 &P_{k2} \ar[r] & K_{k-1,2}   \ar[r] & 0 
}  $$ 
that $F(\iota_{k1})$ is a monomorphism, as desired.  \cqfd

Assume that $F:\A\to \B$ is the zero functor. Then $F$ is totally compatible and the comma category $(F \downarrow \mathcal{B})$ coincides with the product category $\A\times \B$. It follows by Theorem \ref{G-proj/inj} that, for an object $M\in \A\times \B$, $M\in \GP(\A\times \B)$ if and only if $M_1\in \GP(\A)$ and $M_2\cong \Coker\varphi^M\in \GP(\B)$. We record this fact and its dual as a corollary.
\begin{cor} \label{Gproj/inj in AxB} \item 
\begin{enumerate}
    \item If $\A$ and $\B$ have enough projectives, then $\GP(\A\times \B)=\GP(\A)\times \GP(\B)$.
    \item If $\A$ and $\B$ have enough injectives, then $\GI(\A\times \B)=\GI(\A)\times \GI(\B)$.
\end{enumerate}
\end{cor}

Now, we restrict ourselves to triangular matrix rings.

\begin{defn}A $(B,A)$-bimodule $U$ is called:
\begin{enumerate}
    \item left \textbf{totally/weakly compatible} if the functor $U\otimes_A-:A$-Mod $\to B$-Mod is totally/weakly compatible. 
	
	\item left \textbf{totally/weakly cocompatible} if  the functor $\Hom_B(U,-) : B$-Mod $\to A$-Mod is totally/weakly cocompatible.
\end{enumerate}
\end{defn}

In the  following, we provide some examples of such modules.

\begin{prop}\label{tot-comp-exmp} Let $U$ be a $(B,A)$-bimodule. 
	\begin{enumerate}
\item $U$ is left totally compatible if $U_A$ has finite flat dimension and one of the following assertions holds:
\begin{enumerate}
\item $_BU$ has finite projective dimension.
\item $_BU$ has finite FP-injective dimension, and $B$ is left coherent. 
\end{enumerate}
\item  $U$ is left totally cocompatible if $_BU$ has finite injective dimension and one of the following assertions holds:
\begin{enumerate}
\item $U_A$ has finite flat dimension.
\item  $U_A$ has finite FP-injective dimension, and $A$ is right coherent. 
\end{enumerate}
	\end{enumerate}
\end{prop}

\proof (1) Assume that $U_A$ has finite flat dimension.  

Conditions (P1) and (P3) of Definition \ref{P-tot-comp} hold by \cite[Lemma 2.3]{ECT14}.  Note that any totally acyclic $\begin{psmallmatrix}
    \PP_1 \\ \PP_2
\end{psmallmatrix}$ of projective $T$-modules induces an exact complex $\PP_1$ of projective $A$-modules.

To verify condition (P2), let $\QQ$  be a totally acyclic complex of projective $B$-modules and $P$ a projective $A$-module. We claim that $\Hom_B(\QQ,U\otimes_A P)$ is exact.

Since $P$ is a projective $A$-module, it is a direct summand of $A^{(I)}$ for some index set $I$, and so $U\otimes_A P$ is a direct summand of the $B$-module $U\otimes_A A^{(I)} \cong U^{(I)}$.

\textbf{Case (a)} If $_BU$ has projective dimension at most $n$, then so does $U\otimes_A P$, since the class of modules of projective dimension at most $n$ is closed under direct sums and direct summands (Lemma \ref{X_n-closed-dir-sums/summ}). Therefore, $\Hom_B(\QQ,U\otimes_A P)$ is exact, and hence (P2) holds.

\textbf{Case (b)} Assume now that $_BU$ has FP-injective dimension at most $n$, and $B$ is left coherent. It follows that $U\otimes_A P$ has finite FP-injective dimension, as the class of modules with finite FP-injective dimension $\leq n$ is closed under direct sums (since the class of FP-injective modules is) and direct summands (by the dual of Lemma \ref{X_n-closed-dir-sums/summ}). On the other hand, it follows by \cite[Theorem 3]{Gil17} that any FP-injective module belongs to $\GP(B)^\perp$. Since the class $\GP(B)^\perp$ is thick, it contains all modules with finite FP-injective dimension. In particular, it contains $U\otimes_A P$. Since $\QQ$ is a totally acyclic complex of projective $B$-modules, this implies that $\Hom_B(\QQ,U\otimes_A P)$ is exact, as desired.

(2)  The proof is dual. However, since the dual argument requires some additional explanations and references, we include it for completeness.

The conditions (E1) and (E3) of Definition \ref{I-tot-comp} hold by \cite[Lemma 2.4]{ECT14}.  Note that any totally acyclic complex $\begin{psmallmatrix}
		\EE_1\\ \EE_2
		\end{psmallmatrix}$ of injective $T$-modules induces an exact complex $\EE_2$ of injective $B$-modules.

To verify  (E2), let $\EE_1$ be a totally acyclic complex of injective $A$-modules and $E$ an injective $B$-module. We claim that the complex $\Hom_A(\Hom_B(U,E), \EE_1)$ is exact.

Since $_BE$ is injective, it is a direct summand of $(B^+)^{I}$ for some index set $I$, and so  $\Hom_B(U,E)$ is a direct summand of the left $A$-module $$\Hom_B(U,(B^+)^I) \cong \Hom_B(U,(B^+))^I \cong ((B\otimes_B U)^+)^I \cong (U^+)^I.$$ 

\textbf{Case (a)} Assume $U_A$ has flat dimension at most $n$. Then $U^+$ has finite injective dimension $\leq n$ by \cite[Theorem 2.1]{Fie72}. Consequently, $\Hom_B(U,E)$ also has finite injective dimension at most $n$, as the class of modules with injective dimension at most $n$ is closed under direct products and summands (by the dual of Lemma \ref{X_n-closed-dir-sums/summ}).

\textbf{Case (b)} Assume now that $U_A$ has FP-injective dimension at most $n$, and $A$ is right coherent. By \cite[Theorem 2.2]{Fie72}, the left $A$-module $U^+$ has flat dimension at most $n$. Since the class of flat left modules is closed under arbitrary direct products over right coherent rings, $(U^+)^I$ has flat dimension at most $n$, and so does $\Hom_B(U,E)$, as the class of modules with finite flat dimension at most $n$ is closed under direct summands. Now, it follows by \cite[Corollary 5.9]{Sto14} that $^\perp\GI(A)$ contains all modules with finite flat dimension. In particular, it contains $\Hom_B(U,E)$. Since $\EE_1$ is a totally acyclic complex of injective $A$-modules, this implies that $\Hom_A(\Hom_B(U,E), \EE_1)$ is exact, as desired. \cqfd
 
The following result is an immediate consequence of Theorem \ref{G-proj/inj}. It provides a positive answer to Question~1 posed in the Introduction.
\begin{cor} \label{G-proj/inj mod} The following assertions hold.
	\begin{enumerate}
		\item  $\GP(T)=\Phi(\GP(A),\GP(B))$  if and only if $U$ is totally compatible.
	\item $\GI(T)=\Psi(\GI(A),\GI(B))$  if and only if $U$ is totally cocompatible.
	\end{enumerate} 
\end{cor}

We close this section with another application: a characterization of Gorenstein projective/injective objects in the $n$-morphism category $\Mor_n(\C)$.

\begin{cor} \label{G-proj/inj in Mor_n} Let $\C$ be  an abelian category and $n\geq 1$ an integer.

	\begin{enumerate}
		\item If $\C$ has  enough projective objects, then $\GP(\Mor_n(\C))=\Phi_n(\GP(\C))$. 
		\item  If $\C$ has  enough injective objects, then  $\GI(\Mor_n(\C))=\Psi_n(\GI(\C))$.
	\end{enumerate}
		
\end{cor}

\proof (1) By Lemmas \ref{iso of n-Mor} and \ref{iso of Phi_n},  we have isomorphisms of categories
$$\Mor_{n}(\C) \cong (H_n \downarrow \Mor_{n-1}(\C)) \text{ and }\Phi_n(\GP(\C))\cong \Phi(\GP(\C),\Phi_{n-1}(\GP(\C))) \quad \forall n\geq 2$$
with $\Mor_{1}(\C)=\C$ and $\Phi_{1}(\GP(\C))=\GP(\C)$.  The functor $H_n$ is exact and preserves projective objects, so it is totally compatible. Hence, we obtain our desired description of Gorenstein projective objects in $\Mor_n(\C)$, using Theorem \ref{G-proj/inj} and induction. \cqfd

\begin{rem}\label{G-proj/inj in Tn-Mod}  Let $n\geq 2$ be an integer, and let $T_n(R)$ denote the $n \times n$ triangular matrix ring over a ring $R$:
$$T_{n}(R)=\begin{pmatrix}
R      & 0      & \cdots & 0 \\
R      & R      & \ddots & \vdots \\
\vdots & \vdots & \ddots & 0 \\
R      & R      & \cdots & R
\end{pmatrix}$$
 According to \cite[Lemma 1.3.(i)]{Zha11}, the category $T_n(R)$-Mod  of $T_n(R)$-modules is equivalent to the $n$-morphism category $\Mor_n(R\text{-Mod})$.  Therefore, Corollary \ref{G-proj/inj in Mor_n} yields a description of the Gorenstein projective and Gorenstein injective modules over $T_n(R)$ (compare with \cite[Corollaries 4.1 and 4.2]{Zha11}).
    
\end{rem}

\section{Global Gorenstein dimensions of comma categories}\label{Sec. GPD}
 
 In this section we show Theorems B and C from the Introduction. Theorem B will be a particular instance of a more general result that we will prove first.

\bigskip

  \noindent{\bf Setup II.} Besides Setup I from Section \ref{Sec. GP/GI}, we further assume that:
   \begin{enumerate}
      \item[(1)]  $\A$ and $\B$ have direct sums and they are exact whenever we are dealing with the comma category $(F\downarrow \B)$.
      \item[(2)]  $\X_1\subseteq \A$ and $\X_2\subseteq \B$ are classes of objects  projectively resolving  and closed under direct sums, and   $\X:=\Phi(\X_1,\X_2)$  the corresponding class of objects in $(F \downarrow \B)$ is closed under kernels of epimorphisms.
  
      \item[$(1)^{op}$] $\A$ and $\B$ have products, and they are exact whenever we are dealing with the comma category $(\A \uparrow G)$. 
      \item[$(2)^{op}$]  $\Y_1\subseteq \A$ and $\Y_2\subseteq \B$ are classes of objects injectively coresolving and closed under direct products, and  $\Y:=\Psi(\Y_1,\Y_2)$ the corresponding class of objects in  $(\A \uparrow \G)$ is closed under cokernels of monomorphisms. 
  \end{enumerate}
  Note that assumption (2) (resp., $(2)^{op}$), together with \cite[Proposition 2.5]{HZ22}  (resp., its dual) and Lemma \ref{proj-inj in comma} show that $\X$ (resp., $\Y$) is projectively resolving (resp., injectively coresolving). Note also that using Eilenberg’s swindle, it is standard to show  that the classes $\X_1$, $\X_2$, $\Y_1$ and $\Y_2$ are closed under direct summands.

  \bigskip
    \begin{lem}\label{resdim-part-exmp} Let $M_1\in \A$ and $M_2\in \B$ be two objects. Then
 	\begin{enumerate}
\item  $\resdim_{\X_2}(M_2)=\resdim_{\X}\begin{psmallmatrix}
0\\ M_2
\end{psmallmatrix}$.
\item $\resdim_{\X_1}(M_1) \leq \resdim_{\X}\begin{psmallmatrix}
M_1\\ F(M_1)
\end{psmallmatrix}$.
 	\end{enumerate}
 \end{lem}
\proof (1)  Let $n\geq 0$ be an integer and consider a partial projective resolution of $M_2$
$$0\to K_2^n\to P_2^{n-1}\to\cdots\to P_2^0\to M_2 \to 0$$
 Thus, there exists a partial projective resolution of $\begin{psmallmatrix}
0\\M_2
\end{psmallmatrix}$:
$$0\to \begin{psmallmatrix}
0\\K_2^n
\end{psmallmatrix}\to \begin{psmallmatrix}
0\\P_2^{n-1}
\end{psmallmatrix}\to\cdots\to \begin{psmallmatrix}
0\\P_2^0
\end{psmallmatrix}\to \begin{psmallmatrix}
0\\M_2
\end{psmallmatrix}\to 0$$
So we have the following equivalences:
$$  \resdim_{\X}\begin{psmallmatrix}
0\\ M_2
\end{psmallmatrix}  \leq n \Leftrightarrow  \begin{psmallmatrix}
0\\K_2^n
\end{psmallmatrix}  \in  \X   \\
    \Leftrightarrow  K_2^n
 \in  \X_2  \\
    \Leftrightarrow  \resdim_{\X_2}(M_2) \leq n,$$
    which shows $\resdim_{\X}\begin{psmallmatrix}
0\\ M_2
\end{psmallmatrix}=\resdim_{\X_2}(M_2)$.

(2) If $\begin{psmallmatrix}
M_1\\ F(M_1)
\end{psmallmatrix}$ has an $\X$-resolution of length $n$, then $M_1$ has an $\X_1$-resolution of length $n$. This implies that $\resdim_{\X_1}(M_1) \leq \resdim_{\X}\begin{psmallmatrix}
M_1\\ F(M_1)
\end{psmallmatrix}$. \cqfd

 \begin{prop} \label{resdim in comma} Assume $\resdim_{\X_2}(F(\X_1))<\infty$. For any object  $M$ in $(F \downarrow \B)$, we have
 	$$\max\{\resdim_{\X_1}(M_1),\resdim_{\X_2}(M_2)-\resdim_{\X_2}(F(\X_1))\}$$
 	$$\leq \resdim_{\X}(M)
 	\leq$$
 	$$ \max\{\resdim_{\X_1}(M_1)+\resdim_{\X_2}(F(\X_1))+1,\resdim_{\X_2}(M_2)\}$$
 	
 	Consequently, $M$ has finite $\X$-resolution dimension if and only if $M_1$ and $M_2$ have finite   $\X_1$-resolution and  $\X_2$-resolution dimensions, respectively.
 \end{prop}
\proof Let us first prove that if we let $k:=\resdim_{\X_2}(F(\X_1))<\infty$ then $$\max\{\resdim_{\X_1}(M_1),\resdim_{\X_2}(M_2)-k\}\leq \resdim_{\X}(M).$$ We may assume that $n:=\resdim_{\X}(M).$  Then, there exists an exact sequence
  $$0\to X^n\to X^{n-1}\to\cdots\to X^0\to M\to 0$$
  in  $(F \downarrow \B)$ where each $X^i=\begin{psmallmatrix}
  X_1^i\\X_2^i
  \end{psmallmatrix}_{\varphi^{X^i}}\in \X $.
  Thus, there exists an exact sequence 
  $$0\to X_1^n\to X_1^{n-1}\to\cdots\to X_1^0\to M_1\to 0$$
  in $\A$, where each $X_1^i\in \X_1$. So, $\resdim_{\X_1}(M_1)\leq n$. 
  
  For each $i$, there exists an exact sequence in $\B$
  $$0\to F(X_1^i) \stackrel{\varphi^{X^i}}{\longrightarrow}X_2^i\to \Coker\varphi^{X^i}\to 0$$
  where $\Coker\varphi^{X^i}\in \X_2$ . Then, $\resdim_{\X_2}(X^i_2)=\resdim_{\X_2}(F(X_1^i))\leq k$ by Lemma \ref{resdim in ses}(1). So, using the exact sequence 
  $$0\to X_2^n\to X_2^{n-1}\to\cdots\to X_2^0\to M_2\to 0$$ and Lemma  \ref{resdim in ses}(3), we get that 
  $\resdim_{\X_2}(M_2)\leq n+k.$

Now, let us prove the other inequality: $$\resdim_{\X}(M) \leq \max\{\resdim_{\X_1}(M_1)+k+1,\resdim_{\X_2}(M_2)\}.$$ We may assume that $$m:=\max\{\resdim_{\X_1}(M_1)+k+1,\resdim_{\X_2}(M_2)\}< \infty,$$ so if we call $n_1:=\resdim_{\X_1}(M_1)$ and $n_2:=\resdim_{\X_2}(M_2)$ we have $n_1,\ n_2<\infty$. Since $\resdim_{\X_1}(M_1)=n_1\leq m-k-1$, there exists an exact sequence in $\A$ $$0\to X_1^{m-k-1}\to \cdots \to X_1^{1}\stackrel{f_1^1}{\to} X_1^0\stackrel{f_1^0}{\to} M_1\to 0$$
  where each $X_1^i\in \X_1$. Since $\B$ has enough projective objects, there exists an epimorphism $X_2^0\stackrel{g_2^0}{\to} M_2$  with  $X_2^0\in \P(\B)\subseteq \X_2$. Let $K_1^i=\Ker f_1^{i-1}$ and consider the unique morphism   $$f^0_2:F(X_1^0)\oplus X_2^0\to M_2,$$ induced by $\varphi^M F(f^0_1)$ and $g^0_2$. As in the proof of Lemma \ref{enough-proj/inj}, we get a short exact sequence in $(F\downarrow \B)$
    $$0\to K^1\to \pp (X_1^0,X_2^0)\stackrel{f^0}{\to} M\to 0$$
  with $K^1= \begin{psmallmatrix}
  K^1_1\\ K_2^1 \end{psmallmatrix}_{\varphi^{K^1}}$ and $f^0=\begin{psmallmatrix}
  	f_1^0\\f_2^0
  	\end{psmallmatrix}$. Similarly, there exists an epimorphism $X_2^1\stackrel{g_2^1}{\to} K_2^1\to 0$ in $\B$ where $X_2^1\in \P(\B)\subseteq \X_2$, and then 
  we get an exact sequence
  $$0\to K^2\to \pp(X_1^1,X_2^1)\stackrel{}{\to} K^1\to 0.$$
  Repeating this process, we can construct an exact sequence  in $(F\downarrow \B)$
  $$0\to\pp (0,K_2^{m-k}) \to 
   \pp (X_1^{m-k-1},X_2^{m-k-1}) 
  \stackrel{f^{m-k-1}}{\longrightarrow}\cdots  \stackrel{f^{1}}{\longrightarrow} \pp (X_1^0,X_2^0)\stackrel{f^0}{\to} M\to 0.$$

Since each $X_2^i\in \X_2$, it follows by Lemma \ref{X_n-closed-dir-sums/summ} that 
    $$\resdim_{\X_2}(F(X_1^i)\oplus X_2^i)=\resdim_{\X_2}( F(X_1^i))\leq k.$$
So, applying Lemma \ref{resdim in ses}(2) to the exact sequence
 $$0\to K^{m-k}_2
  \to F(X_1^{m-k-1})\oplus X_2^{m-k-1}
  \stackrel{f_2^{m-k-1}
  }{\longrightarrow}\cdots\to F(X_1^0)\oplus X_2^0
  \stackrel{f_2^0
  }{\to} M_2\to 0,$$ we get that $\resdim_{\X_2}(K_2^{m-k})\leq \max\{k,n_2-m+k\}\leq k$ since $n_2\leq m$. This means that there exists an exact sequence 
  $$0\to X_2^m\to \cdots\to X_2^{m-k+1}\to X_2^{m-k}\to K_2^{m-k}\to 0$$
 in $\B$ with each $X_2^i\in \X_2.$ Thus, there exists an exact sequence in $(F \downarrow \B)$ $$
\begin{tikzcd}[column sep=1em, row sep=3em]
  0 \arrow[r] & \pp(0,X_2^{m}) \arrow[r, "f^m"]
  &  \cdots \arrow[r]
  &  \pp(0,X_2^{m-k+1}) \arrow[rr, "f^{m-k+1}"]  & & \pp(0,X_2^{m-k}) \arrow[dllll, out=-90, in=150, looseness=0.45] & 
  \\
  & \cdots \arrow[r]
  & \pp(X_1^1,X_2^1) \arrow[r, "f^1"]
  & \pp (X_1^0,X_2^0) \arrow[r, "f^0"]
  & M \arrow[r]
  & 0. 
\end{tikzcd}
$$
Clearly, all $\pp(X_1^i,X_2^i),\pp(0,X_2^j)\in \X$.  Thus, $\resdim_\X(M)\leq m$. \cqfd

\begin{thm} \label{glresdim of comma}

$$\max\{\glresdim_{\X_1}(\A),\glresdim_{\X_2}(\B)\}$$
$$\leq \glresdim_{\X}(F \downarrow \B)
\leq$$
$$ \max\{\glresdim_{\X_1}(\A)+\resdim_{\X_2}(F(\X_1))+1,\glresdim_{\X_2}(\B)\}$$
 \end{thm}

\proof We only show the inequality $$\max\{\glresdim_{\X_1}(\A),\glresdim_{\X_2}(\B)\}\leq \glresdim_{\X}(F \downarrow \B),$$ 
as the other inequality follows from Proposition \ref{resdim in comma}. 

We may assume that $n:=\glresdim_{\X}(F \downarrow \B)<\infty.$ For any objects  $M_1\in \A$ and $M_2\in \B$, we have by Lemma \ref{resdim-part-exmp} 
$$\resdim_{\X_2}(M_2)=\resdim_\X\begin{psmallmatrix}
0\\M_2
\end{psmallmatrix}\leq n \text{ and } \resdim_{\X_1}(M_1)\leq \resdim_\X\begin{psmallmatrix}
M_1\\F(M_1)
\end{psmallmatrix}\leq n$$
Thus, $\glresdim_{\X_1}(\A)\leq n$ and  $\glresdim_{\X_2}(\B)\leq n$.  \cqfd

 Dually, we have the corresponding dual results of Proposition \ref{resdim in comma} and Theorem \ref{glresdim of comma}.
 
 \begin{prop} \label{coresdim in comma} Assume $\coresdim_{\Y_1}(G(\Y_2))<\infty$. For any object $M$ in $(\A \uparrow \G)$, we have
    $$\max\{\coresdim_{\Y_1}(M_1) -\coresdim_{\Y_1}(G(\Y_2)),\coresdim_{\Y_2}(M_2)\}$$ $$ \leq  \coresdim_\Y(M) \leq $$
    $$\max\{\coresdim_{\Y_1}(M_1),\coresdim_{\Y_2}(M_2)+ \coresdim_{\Y_1}(G(\Y_2))+1\}.$$
 \end{prop}
 \begin{thm} \label{glcoresdim of comma}
 	
 	$$\max\{\glcoresdim_{\Y_1}(\A),\glcoresdim_{\Y_2}(\B)\}$$ $$ \leq  \glcoresdim_\Y(\A \uparrow \G) \leq $$
    $$\max\{\glcoresdim_{\Y_1}(\A),\glcoresdim_{\Y_2}(\B) + \coresdim_{\Y_1}(G(\Y_2))+1\}$$
 \end{thm}

\bigskip

The homological dimensions $\resdim_{\X_2}(F(\X_1))$ and $\coresdim_{\Y_1}(G(\Y_2))$ play a central role in the estimates established in Propositions~\ref{resdim in comma} and \ref{coresdim in comma} and Theorems~\ref{glresdim of comma} and \ref{glcoresdim of comma}. In the Gorenstein setting, these quantities deserve special attention, as they admit a natural interpretation: they measure how far a functor is from preserving Gorenstein objects.

\begin{defn}
Let $F:\A\to\B$ and $G:\B\to\A$ be functors, $A$ and $B$ two rings, and $U$ be a $(B,A)$-bimodule.

\begin{enumerate}
\item The \textbf{global Gorenstein projective dimension} of $F$ is defined as:
\[
\GPD(F):=\sup\{\Gpd_\B(F(A))\mid A\in\GP(\A)\}.
\]
In particular, the \textbf{left global Gorenstein projective dimension} of $U$ is defined as:
\[
\lGPD(U):=\GPD(U\otimes_A-).
\]

\item The \textbf{global Gorenstein injective dimension} of $G$ is defined as:
\[
\GID(G):=\sup\{\Gid_\A(G(B))\mid B\in\GI(\B)\}.
\]
In particular, the \textbf{left global Gorenstein injective dimension} of $U$ is defined as:
\[
\lGID(U):=\GID(\Hom_B(U,-)).
\]
\end{enumerate}
\end{defn}

The global Gorenstein projective (resp., injective) dimension of an abelian category measures how far the category is from having  the property that every object is Gorenstein projective (resp., injective). Analogously, the global Gorenstein projective (resp., injective) dimension of a functor measures how far the functor is from preserving Gorenstein projective (resp., injective) objects:
\begin{enumerate}
    \item $F$ preserves Gorenstein projective objects if and only if  $\GPD(F)=0$. 
     \item $G$ preserves Gorenstein injective objects if and only if  $\GID(G)=0$.
     \end{enumerate}

\begin{rem}One could similarly define the global projective/injective dimensions of functors:
     $$\PD(F):=\sup\{ \pd_\B(F(A)) \,|\, A\in \P(\A) \},$$
     $$\ID(G):=\sup\{ \id_\A(G(B)) \,|\, B\in \I(\B) \}.$$

Under some assumptions, the global projective/injective dimensions of functors are nothing but the projective/injective dimensions of some objects:
\begin{enumerate}
\item  If $\A$ has a projective generator $P$  and $F$ preserves  direct sums, then we have $\PD(F)=\pd_\B(F(P))$.
      
\item If $\B$ has an injective cogenerator $E$ and $G$ preserves direct products, then $\ID(G)=\id_\A(G(E)).$
\end{enumerate}
Assertion (2) is dual. To see (1), we only need to argue the inequality $\PD(F)\leq \pd_\B(F(P))$,  since $P$ is projective.

Let $Q$ be a projective object of $\A$. Since $P$ is a projective generator,
there exists an index set $I$ such that
$Q$ is a direct summand of $P^{(I)}$.
Because $F$ preserves coproducts, $
F(P^{(I)})\cong F(P)^{(I)}.$ Hence, $F(Q)$ is a direct summand of $F(P)^{(I)}$, and therefore
$\pd_\B(F(Q))
\leq
\pd_\B(F(P)^{(I)})
=
\pd_\B(F(P)).
$ Since this holds for every projective object $Q$, taking the supremum yields $
\PD(F)\leq\pd_\B(F(P)),$ which completes the proof. 
\end{rem}

The class of Gorenstein projective objects in an abelian category with enough projective (resp., injective) objects is  projectively resolving (resp., injectively coresolving), closed under direct sums (resp., products) and direct summands; see e.g., \cite[Theorem 2.5]{Hol04} for module categories or \cite[Theorem 11.2.4]{Per16} for abelian categories. Using  Theorem~\ref{G-proj/inj},  the following  estimates of the Gorenstein dimension of objects in comma categories are an immediate application of Propositions~\ref{resdim in comma} and \ref{coresdim in comma} to the classes of Gorenstein projective and injective objects, respectively.
 \begin{prop}  Let $M\in (F \downarrow \B)$ and $N\in (\A \uparrow \G)$ be two objects.  
\begin{enumerate}
\item Assume  $F:\A \to \B$ is totally compatible with finite global Gorenstein  projective dimension. Then
$$\max\{\Gpd_\A(M_1),\Gpd_\B(M_2)-\GPD(F)\} $$ $$\leq \Gpd_{(F \downarrow \B)}(M)
\leq $$  $$\max\{\Gpd_\A(M_1)+\GPD(F)+1,\Gpd_\B(M_2)\}$$
Consequently, $M$ is of finite Gorenstein  projective dimension if and only if so are $M_1$ and $M_2$.
\item Assume $G:\B\to \A$ is totally cocompatible with finite global Gorenstein injective dimension. Then 
$$\max\{\Gid_\A(N_1)-\GID(G),\Gid_\B(N_2)\}$$
$$\leq \Gid_{(\A \uparrow \G)}(N)
\leq $$
 $$\max\{\Gid_\A(N_1),\Gid_\B(N_2)+\GID(G)+1\}$$
Consequently, $N$ is of finite  Gorenstein injective dimension if and only if so are $N_1$ and $N_2$.

\end{enumerate} 

\end{prop}
Now we are ready to prove Theorem B from the Introduction. Using Theorem \ref{G-proj/inj} again, the proof is an immediate application of Theorems \ref{glresdim of comma} and \ref{glcoresdim of comma} to the classes of Gorenstein projective and injective objects, respectively.
\begin{thm}  \label{global-estimation} \item   
\begin{enumerate}
\item Assume  $F:\A \to \B$ is totally compatible. Then
$$\max\{\GPD(\A),\GPD(\B)\}\leq \GPD(F \downarrow \B)
\leq \max\{\GPD(\A)+\GPD(F)+1,\GPD(\B)\}$$
Consequently, the comma category  $(F \downarrow \B)$ is of finite global Gorenstein  projective dimension if and only if so are $\A$ and $\B$.
\item Assume $G:\B\to \A$ is totally cocompatible. Then 
$$\max\{\GID(\A),\GID(\B)\}
\leq \GID(\A \uparrow \G)
\leq
 \max\{\GID(\A),\GID(\B) +\GID(G)+1\}$$
Consequently, the comma category $(\A \uparrow \G)$  is of finite global Gorenstein injective dimension if and only if so are $\A$ and $\B$.

\end{enumerate} 

\end{thm}

\begin{rem} Keeping in mind Lemma \ref{proj-inj in comma}, Theorems \ref{glresdim of comma} and \ref{glcoresdim of comma} can also be applied to the classes of projective and injective objects, respectively:
$$\max\{\PD(\A),\PD(\B)\}\leq \PD(F \downarrow \B)
\leq \max\{\PD(\A)+\PD(F)+1,\PD(\B)\}$$
$$\max\{\ID(\A),\ID(\B)\}
\leq \ID(\A \uparrow \G)
\leq
 \max\{\ID(\A),\ID(\B) +\ID(G)+1\}$$
where, for an abelian category $\C$, $\PD(\C):=\sup\{\pd_\C(M)\mid M\in\C \}$ and $\ID(\C):=\sup\{\id_\C(M)\mid M\in\C \}$ denote the global projective and the global injective dimension of $\C$, respectively. 

As a consequence, the comma category $(F \downarrow \B)$ (resp., $(\A \uparrow \G)$) is of finite global projective (resp., injective) dimension if and only if so are $\A$ and $\B$.
\end{rem}

Now we prove Theorem C from the Introduction, that is, we provide several situations in which the estimates in Theorem \ref{global-estimation}  become exact formulas, giving a partial answer to Question 2 raised in the Introduction. 

	\begin{prop} \label{Formula 1}\item 
    \begin{enumerate}
        \item If $F=0$, then $\GPD(F \downarrow \B)=\GPD(\A\times \B)= \max\{\GPD (\A),\GPD(\B)\}.$
        \item If $G=0$, then $\GID(\A \uparrow \G)=\GID(\A\times \B)= \max\{\GID (\A),\GID(\B)\}.$
    \end{enumerate}
\end{prop}

\proof  Assume $F=0$. In this case,  $(F \downarrow \B)$ is nothing but the product category $\A\times \B$. By Corollary \ref{Gproj/inj in AxB}, for an object  $(M_1,M_2)\in \A\times \B$, $(M_1,M_2)\in \GP(\A\times \B)$ if and only if $M_1\in \GP(\A)$ and $M_2\in \GP(\B)$. This leads immediately to 
$$\Gpd_{\A\times \B}(M_1,M_2)=\sup\{\Gpd_\A(M_1),\Gpd_\B(M_2)\}.$$
For any integer $n \geq 0$, we have the equivalences:
\begin{eqnarray*}
   \GPD(\A\times \B) \leq n &\Leftrightarrow& \forall (M_1,M_2)\in \A\times \B, \,  \Gpd_{\A\times \B}(M_1,M_2)\leq n \\
   &\Leftrightarrow& \forall M_1\in \A,  \,  M_2\in \B, \Gpd_{\A}(M_1)\leq n \text{ and } \Gpd_{\B}(M_2)\leq n \\
 &\Leftrightarrow&   \GPD(\A) \leq n \text{ and } \GPD(\B) \leq n.
\end{eqnarray*}
This shows our desired equality. \cqfd

Recall that an abelian category is called Frobenius if it has enough projective objects and enough injective objects, and the class of projective objects coincides with the class of injective objects. The following lemma is an immediate consequence of Lemma \ref{glGdim-charac}.

	\begin{lem} \label{Frob-categ-lem} Assume that  $\mathcal{C}$ is an abelian category with enough projectives and injectives. Then $\C$ is  Frobenius if and only if $\glGdim(\C)=0$.
\end{lem}

\begin{cor} Assume that $\A$ and $\B$ have both enough projective and injective objects. 

\item \begin{enumerate}
    \item If $F$  satisfies (P3), then  $(F \downarrow \B)$ is Frobenius if and only if so are $\A$ and $\B$ and $F=0$.
    \item  If $G$ satisfies (E3), then  $(\A \uparrow \G)$ is Frobenius if and only if so are $\A$ and $\B$ and $G=0$.
\end{enumerate}

\end{cor}

\proof By Lemma \ref{Frob-categ-lem}, we only need to show that $\GPD(F \downarrow \B)=0$ if and only if $\GPD(\A)=\GPD(\B)=0$ and $F=0$. However, this will follow from Proposition \ref{Formula 1} once we show that $F=0$ whenever $\GPD(F \downarrow \B)=0$.

Assume $\GPD(F \downarrow \B)=0$. By Lemma \ref{Frob-categ-lem}, every object in $(F \downarrow \B)$ is Gorenstein projective. This implies, in particular, that $\pp$ preserves Gorenstein projective objects, and hence $F$ is weakly compatible by Lemma \ref{Goren-preserv}. Since it satisfies (P3), it is totally compatible.

Assume now that $F\neq 0$. So we get an object $M_1\in\A$ such that $F(M_1)\neq 0$. It follows by Theorem \ref{G-proj/inj} that the object $M:=\begin{psmallmatrix}
    M_1\\ 0
\end{psmallmatrix}$ is not Gorenstein projective as $\Ker\varphi^M=F(M_1) \neq 0$, which contradicts the fact that every object in $(F \downarrow \B)$ is Gorenstein projective. Hence $F=0$. \cqfd

Recall that a ring $R$ is called quasi-Frobenius if injective left (resp.,  right) $R$-modules coincide with projective left (resp., right) $R$-modules. Equivalently, $R$ is quasi-Frobenius if and only if $R$-Mod is a Frobenius category. 

\begin{cor} Assume $U\otimes_A-$  satisfies (P3) or $\Hom_B(U,-)$ satisfies (E3). Then $T$ is quasi-Frobenius if and only if so are $A$ and $B$ and $U=0$.
    
\end{cor}

	\begin{prop} \label{Formula 2} \item 
    \begin{enumerate}
	    \item  If $\A$ is Frobenius and $0\neq F:\A \to \B$ is weakly compatible, then 
           $$\GPD(F \downarrow \B)=
			\max\{\GPD(F)+1,\GPD(\B)\}.$$

            \item  If $\B$ is Frobenius and $0\neq G:\B \to \A$ is weakly cocompatible, then  $$ \GID(\A \uparrow G)= \max\{\GID(\A),\GID(G)+1\}.$$
        
	\end{enumerate}
	\end{prop} 
\proof Set $\T=(F \downarrow \B)$  and assume that $\A$ is Frobenius. First, note that $F:\A \to \B$ is totally compatible (see Remark \ref{rem:tot-comp=comp}). By Lemma \ref{Frob-categ-lem}, we obtain $\GPD(\A)=0$, and so it follows by Theorem \ref{global-estimation} that 
$$\GPD(\B)\leq \GPD(F \downarrow \B)
\leq \max\{\GPD(F)+1,\GPD(\B)\}.$$
If $\GPD(F)+1\leq \GPD(\B)$, then $\GPD(F \downarrow \B)=\GPD(\B)$, and we are done. So, let us assume $\GPD(\B) < \GPD(F)+1$. As $\GPD(F)\leq \GPD(\B)$, we obtain that $n=\GPD(F)= \GPD(\B).$ 

Let $M_1\in \GP(\A)=\A$ be an object such that $n=\Gpd_\B(F(M_1))$. If $n=0$, we may choose   $M_1$ with $F(M_1) \neq  0$, which exists because $F \neq  0$. Then $M:=\begin{psmallmatrix}
    M_1 \\ 0
\end{psmallmatrix}$ is not Gorenstein projective  because $\varphi^M$ is not a monomorphism, so $\Gpd_\T M=1= n+1$. 

Assume $n\geq 1$ and  consider the short exact sequence in $(F \downarrow \B)$
$$0\to \begin{psmallmatrix}
     0 \\ F(M_1)
\end{psmallmatrix}\to \begin{psmallmatrix}
   M_1 \\ F(M_1)
\end{psmallmatrix}\to \begin{psmallmatrix}
    M_1 \\ 0
\end{psmallmatrix}\to 0.$$
By Theorem \ref{G-proj/inj},  $\begin{psmallmatrix}
   M_1 \\ F(M_1)
\end{psmallmatrix}$ is Gorenstein projective. It follows from Lemmas \ref{resdim in ses}(3) and  \ref{resdim-part-exmp}(1) that 
$$\Gpd_{(F \downarrow \B)}\begin{psmallmatrix}
    M_1 \\ 0
\end{psmallmatrix}=\Gpd_{(F \downarrow \B)}\begin{psmallmatrix}
     0 \\ F(M_1)
\end{psmallmatrix}+1=\Gpd_\B(F(M_1))+1=n+1.$$
This shows that $\GPD(F \downarrow \B)=n+1=\GPD(F)+1$. \cqfd

\begin{cor} Assume $U\neq 0$.
\begin{enumerate}
    \item If $A$ is quasi-Frobenius and $U$ is weakly compatible, then 
  $$\glGdim(T)=\max\{\lGPD(U)+1,\glGdim(B)\}.$$
  \item If $B$ is quasi-Frobenius and $U$ is weakly cocompatible, then 
  $$\glGdim(T)= \max\{\glGdim(A),\lGID(U)+1\}.$$
\end{enumerate}
\end{cor}

Frobenius pairs and functors are very classical in algebra.  Let $\alpha:\A \to \A$ and $\beta:\B \to \B$ be two auto-equivalences. Recall that an adjoint pair $(F,G)$ is called  $(\alpha,\beta^{-1})$-strongly Frobenius in the sense of \cite{Mor65} or a Frobenius pair of type $(\alpha,\beta)$ in the sense of \cite{CR22} if both $(F, G)$ and $(G, \beta F \alpha)$ are adjoint pairs. Here we simply call such a pair a Frobenius pair. One of the main results of Chen and Ren in \cite{CR22} states that faithful Frobenius functors between abelian categories preserve and reflect objects with bounded Gorenstein projective dimension. Precisely, if $F:\A \to \B$ is a faithful Frobenius functor, then we have $\Gpd_\A(X)=\Gpd_\B(F(X))$ for each object $X\in \A$; see \cite[Theorem 3.2]{CR22}. As an immediate consequence, we have $\glGdim(\A)\leq \glGdim(\B)$. This result will be key to obtaining our next formula.

 \begin{lem}\label{Frob-vs-comp} Assume $(F,G)$ is a Frobenius pair. Then 
 \begin{enumerate}
     \item $F$ is totally compatible and $\GPD(F)=0$.
     \item $G$ is totally cocompatible and $\GID(G)=0$.
 \end{enumerate}
 \end{lem}
\proof (1) The functor $F$ is totally compatible, as it is exact and preserves projective objects by \cite[Corollary 2.2]{CR22}. As $F$ also preserves Gorenstein projective objects by \cite[Lemma 3.6]{CR22}, $\GPD(F)=0$. 

(2) is dual. \cqfd

  Recall that if $(F,G)$ is an adjoint pair, then $(F \downarrow \B)$ and $(\A \uparrow \G)$ are isomorphic, and both are denoted by $\T$.

\begin{prop}  \label{Formula 3} Assume $(F,G)$ is a Frobenius pair such that $F$ and $G$ are nonzero functors.
\begin{enumerate}
    \item If   $ F : \A \to \B$ is  faithful, then
        $ \GPD(\T)=\max\{\GPD(\A)+1,\GPD(\B)\}.$
    \item If   $ G : \B \to \A$ is  faithful, then
        $\GID(\T)=\max\{\GID(\A),\GID(\B)+1\}.$
\end{enumerate}

\end{prop}
\proof By Lemma \ref{Frob-vs-comp}, $F$ is totally compatible and  $\GPD(F)=0$. So, it follows by Theorem \ref{global-estimation} that we have 
$$\max\{\GPD(\A),\GPD(\B)\}\leq \GPD(\T) \leq  \max\{\GPD(\A)+1,\GPD(\B)\}.$$
Using \cite[Theorem 3.2]{CR22} again, we obtain that $\GPD(\A)\leq \GPD(\B)$. Clearly $\GPD(\T)=\GPD(\B)$ when $\GPD(\A)<\GPD(\B)$. So, we may assume that $n:=\GPD(\A)=\GPD(\B)<\infty$. Then $ n\leq \GPD(\T) \leq n+1.$

Let $M_1\in \A$ be an object such that $n=\Gpd_\A(M_1)$. If $n=0$, we may choose   $M_1$ with $F(M_1) \neq  0$, which exists because $F \neq  0$. Then $M:=\begin{psmallmatrix}
    M_1 \\ 0
\end{psmallmatrix}$ is not Gorenstein projective  because $\varphi^M$ is not a monomorphism, so $\Gpd_\T M=1= n+1$.

Assume $n\geq 1$. Using \cite[Theorem 3.2]{CR22}, we obtain that  $\Gpd_\B(F(M_1))=\Gpd_\A(M_1)=n$. Then  $\Ext_\B^n(F(M_1),Q)\neq 0$ for some projective object $Q\in \B$ (see, e.g.,  \cite[Theorem 2.20]{Hol04}). Applying the functor $\Hom_{\T}(-,\begin{psmallmatrix}
0\\Q
\end{psmallmatrix})$ to the short exact sequence
$$0\to \begin{psmallmatrix}
     0 \\ F(M_1)
\end{psmallmatrix}\to \begin{psmallmatrix}
   M_1 \\ F(M_1)
\end{psmallmatrix}\to \begin{psmallmatrix}
    M_1 \\ 0
\end{psmallmatrix}\to 0,$$
we get an exact sequence of abelian groups
$$\begin{tikzcd}[column sep=1em, row sep=3em]
  \cdots  \arrow[r] & \Ext^n_{\T}(\begin{psmallmatrix}
M_1\\F(M_1)
\end{psmallmatrix},\begin{psmallmatrix}
0\\Q
\end{psmallmatrix}) \arrow[r]
  &  \Ext^n_{\T}(\begin{psmallmatrix}
0\\F(M_1)
\end{psmallmatrix},\begin{psmallmatrix}
0\\Q
\end{psmallmatrix}) \arrow[dl, out=-16, in=150, looseness=0.8] 
  \\
    &   \Ext^{n+1}_{\T}(\begin{psmallmatrix}
M_1\\0
\end{psmallmatrix},\begin{psmallmatrix}
0\\Q
\end{psmallmatrix}) \arrow[r] & \Ext^{n+1}_{\T}(\begin{psmallmatrix}
M_1\\F(M_1)
\end{psmallmatrix},\begin{psmallmatrix}
0\\Q
\end{psmallmatrix}) \arrow[r]
  & \cdots 
\end{tikzcd}
$$
By Lemma \ref{Ext-iso}, we obtain for each $i\geq 1$, $\Ext^{i}_{\T}(\begin{psmallmatrix}
M_1\\F(M_1)
\end{psmallmatrix},\begin{psmallmatrix}
0\\Q
\end{psmallmatrix})=0$. So using Lemma \ref{Ext-iso} again and the above long exact sequence, we obtain the isomorphisms $$\Ext^{n+1}_{\T}(\begin{psmallmatrix}
M_1\\0
\end{psmallmatrix},\begin{psmallmatrix}
0\\Q
\end{psmallmatrix})\cong \Ext^n_{\T}(\begin{psmallmatrix}
0\\F(M_1)
\end{psmallmatrix},\begin{psmallmatrix}
0\\Q
\end{psmallmatrix})\cong \Ext_\B^n(F(M_1),Q)\neq 0.$$ 
Since $\begin{psmallmatrix}
0\\Q
\end{psmallmatrix}\in \T$ is projective, $\Gpd_\T\begin{psmallmatrix}
M_1\\0
\end{psmallmatrix}\geq n+1$.  But $\Gpd_\T\begin{psmallmatrix}
M_1\\0
\end{psmallmatrix}\leq \GPD(\T)
\leq n+1$, so it follows that $n<\Gpd_\T\begin{psmallmatrix}
M_1\\0
\end{psmallmatrix}$. Thus $\Gpd_\T\begin{psmallmatrix}
M_1\\0
\end{psmallmatrix}=n+1$. This shows that $\GPD(\T)
= n+1= \GPD(\A)+1.$
\cqfd

We close our paper with one last result relating the global Gorenstein dimension of the $n$-morphism category with that of the base category.

	\begin{prop}  \label{Formula 4: Mor_n} Let $\C$ be an abelian category and $n\geq 2$ an integer.
    \begin{enumerate}
        \item If $\C$ has enough projectives, and direct sums and they are exact, then
        $$\GPD(\Mor_{n}(\C))
		=\GPD(\C)+1.$$
         \item If $\C$ has enough injectives, and direct products and they  are exact, then
        $$\GID(\Mor_{n}(\C))
		=\GID(\C)+1.$$
    \end{enumerate}
In particular, if $\C$ satisfies the assumptions given (1) and (2), then
        $$\glGdim(\Mor_{n}(\C))
		=\glGdim(\C)+1.$$
		
	\end{prop}

    \proof  We proceed by induction on $n$. To do so, we use  the isomorphism of categories from Lemma \ref{iso of n-Mor} $$\Mor_{n}(\C) \cong (H_n \downarrow \Mor_{n-1}(\C)).$$ 

   Let $n=2$ and let us prove $\GPD(\Mor_2(\C))
		=\GPD(\C)+1.$ (Note that $\Mor_{1}(\C)=\C$).  In this case, the functor $H_2$ is the identity $\id_\C:\C\to \C$, so it is a faithful Frobenius functor. By Proposition \ref{Formula 3}, we have 
$$\GPD(\Mor_2(\C))=\max\{\GPD(\C)+1,\GPD(\C)   \}=\GPD(\C)+1.$$

 Assume now $n\geq 3$. As the functor $H_n$ preserves Gorenstein projective objects, by Corollary \ref{G-proj/inj in Mor_n}, it follows by Theorem \ref{global-estimation}(1)
$$\max\{\GPD(\C),\GPD(\Mor_{n-1}(\C))\}$$
$$\leq \GPD(\Mor_{n}(\C))
\leq $$
$$ \max\{\GPD(\C)+1,\GPD(\Mor_{n-1}(\C))\}$$
 By the induction hypothesis, we have $\GPD(\Mor_{n-1}(\C))
		=\GPD(\C)+1$. Hence, 
$$ \GPD(\Mor_{n}(\C))=\GPD(\Mor_{n-1}(\C))=\GPD(\C)+1.$$
        \cqfd

        \begin{rem}
            Keeping in mind Remark \ref{G-proj/inj in Tn-Mod}, we have, by Proposition \ref{Formula 4: Mor_n}, the following formula:
         $$\glGdim(T_{n}(R))
		=\glGdim(R)+1.$$
        \end{rem}

Finally, we conclude the paper with a natural question that the reader may ask: what can be said about Gorenstein flat modules and dimensions? This question is the subject of another work currently in preparation \cite{EGO26}, in which we investigate rings satisfying a Gorenstein analogue of Lambek's theorem: a module $M$ is Gorenstein flat if and only if its character module $M^+:=\Hom_{\mathbb{Z}}(M,\mathbb{Q}/\mathbb{Z})$ is Gorenstein injective. In the course of that investigation, we were led to study Gorenstein flat modules and dimensions over triangular matrix rings. Thus, the present paper provides the foundation for that work, while the two papers complement each other.

\bigskip

\noindent {\bf Acknowledgment.} The authors were partially supported by a project from the Spanish Ministerio de Ciencia, Innovación y Universidades through its Agencia Estatal de Investigación, REFERENCIA DEL PROYECTO/AEI/PID2024-158993NB-I00.


\bigskip




	\begin{thebibliography}{999}

\bibitem[ARS95]{ARS95}
M. Auslander, I. Reiten and S.O. Smal\o, \textit{Representation theory of Artin algebras}, Cambridge Stud. Adv. Math. \textbf{36}, Cambridge Univ. Press, 1995.
	
\bibitem[BEGO21]{BEGO21}
D. Bennis, R. El Maaouy,  J. R. Garc\'{\i}a Rozas and L. Oyonarte,  \textit{Relative Gorenstein dimensions over triangular matrix rings},  Mathematics \textbf{9} (2021), 2676.

\bibitem[BM10]{BM10}
D. Bennis and  N. Mahdou, \textit{Global Gorenstein dimensions}, Proc. Am. Math. Soc. \textbf{138} (2010), 461--465. 


 \bibitem[CR22]{CR22}
    X.W. Chen and W. Ren, \textit{Frobenius functors and Gorenstein homological properties}, J. Algebra \textbf{610} (2022), 18--37.


 \bibitem[Elm26]{Elm26} 
 R. El Maaouy, \textit{Homological algebra in categories of representations of quivers via duality pairs}, J. Algebra \textbf{708} (2026), 91--137.

\bibitem[EGO26]{EGO26}
 R. El Maaouy,  J. R. Garc\'{\i}a Rozas and L. Oyonarte,  \textit{Rings over which a module is Gorenstein flat if and only if its character module is Gorenstein injective}, in preparation.

	\bibitem[ECT14]{ECT14}
	E. Enochs, M. Cort\'es-Izurdiaga and B. Torrecillas, \textit{Gorenstein conditions over triangular matrix rings}, J. Pure Appl. Algebra \textbf{218} (2014), 1544--1554.

	
	\bibitem[EJ00]{EJ00} 
	E. E. Enochs and O. M. G. Jenda, \textit{Relative Homological Algebra}, de Gruyter Expositions in Mathematics. Vol. \textbf{30}, Walter de Gruyter (2000).
	


    \bibitem[EHSL16]{EHSL16}
    H. Eshraghi, R. Hafezi, Sh. Salarian and  Z.W. Li, \textit{Gorenstein projective modules over triangular matrix rings}, Algebra Colloq.  \textbf{23} (2016), 97--104.

	
	\bibitem[Fie72]{Fie72}
	D. J. Fieldhouse, \textit{Character modules, dimension and purity}, Glasgow Math. J. \textbf{13} (1972), 144--146.
	
	\bibitem[FGR75]{FGR75}
	R.M. Fossum, P. Griffith and I. Reiten, \textit{Trivial extensions of abelian categories. Homological algebra of trivial extensions of abelian categories with applications to ring theory}, Lecture Notes in Mathematics 456, Springer-Verlag, 1975.

    \bibitem[Gil17]{Gil17}
    J. Gillespie, \textit{On Ding injective, Ding projective and Ding flat modules and complexes}, Rocky Mountain J. Math. \textbf{47} (2017), 2641--2673.
	
	
	
	\bibitem[Gre82]{Gre82}
	E.L. Green, \textit{On the representation theory of rings in matrix form}, Pacific J. Math. \textbf{100} (1982), 123--138.


	
	\bibitem[HV99]{HV99}
	A. Haghany and K. Varadarajan, \textit{Study of formal triangular matrix rings}, Comm. Algebra \textbf{27} (1999), 5507--5525.
	
	\bibitem[HV00]{HV00}
	A. Haghany and K. Varadarajan, \textit{Study of modules over formal triangular matrix rings}, J. Pure Appl. Algebra \textbf{147} (2000), 41--58.
	
	\bibitem[Hol04]{Hol04} 
	H. Holm, \textit{Gorenstein homological dimensions}, J. Pure Appl. Algebra \textbf{189} (2004), 167--193.
	
	
	\bibitem[HZ22]{HZ22}
	J. Hu and H. Zhu, \textit{Special precovering classes in comma categories}, Sci. China Math.   \textbf{65} (2022), 933--950.
	
	\bibitem[Hua22]{Hua22}
	Z. Huang, \textit{Homological dimensions relative to preresolving subcategories II}, Forum Math. \textbf{34} (2022), 507--530.

    \bibitem[LZHZ20]{LZHZ20}
    H. H. Li, Y. F. Zheng, J. S. Hu, and H. Y. Zhu, \textit{Gorenstein projective modules and recollements over triangular matrix rings}, Comm. Algebra \textbf{48} (11) (2020), 4932--4947.

    	\bibitem[Mao20]{Mao20}
L. Mao, \textit{Cotorsion pairs and approximation classes over formal triangular matrix rings}, J. Pure Appl. Algebra \textbf{224}  (2020), 106271.


   \bibitem[Mao22]{Mao22}
    L. Mao, \textit{Ding modules and dimensions over formal triangular matrix rings}, Rend. Semin. Mat. Univ. Padova \textbf{148} (2022), 1--22.

    \bibitem[MY20]{MY20}
    H. Minamoto and K. Yamaura, \textit{Homological dimension formulas for trivial extension algebras}, J. Pure Appl. Algebra \textbf{224} (2020), 106344.
    

\bibitem[Mor65]{Mor65}
K. Morita, \textit{Adjoint pairs of functors and Frobenius extensions}, Sci. Rep. Tokyo Kyoiku Daigaku, Sect. A \textbf{9} (1965), 40--71.
    



   \bibitem[PR73]{PR73}
    I. Palmer and J.-E. Roos, \textit{Explicit formulae for the global homological dimensions of trivial extensions of rings}, J. Algebra \textbf{27} (1973), 380--413.


    \bibitem[PZH22]{PZH22} 
    Y. Y. Peng, R. M. Zhu and Z. Y. Huang, \textit{Gorenstein projective objects in comma categories}, Period. Math. Hung. \textbf{84} (2022), 186--202.

    \bibitem[P\'er16]{Per16} 
      M. A. P\'erez, \textit{Introduction to Abelian Model Structures and Gorenstein Homological Dimensions}. Monographs and Research Notes in Mathematics. CRC Press, Boca Raton, FL. (2016)

      \bibitem[Rot79]{Rot79} 
    J. J. Rotman, \textit{An Introduction to Homological Algebra} (Academic Press, New York, 1979).

\bibitem[\v{S}to14]{Sto14}
J. \v{S}\v{t}ov\'{\i}\v{c}ek, \textit{On purity and applications to coderived and singularity categories}, arXiv:1412.1615.
	

  \bibitem[WL22]{WL22}
    Z. P. Wang and Z.K. Liu, \textit{Recollements induced by monomorphism categories}, J. Algebra \textbf{594} (2022), 614--635.

     \bibitem[YHW24]{YHW24} 
    Y. Yuan, J. He and D. Wu, \textit{Cotorsion pairs in comma categories}, Czechoslovak Math. J. \textbf{74} (149) (2024), 715--734.
    
\bibitem[Zha11]{Zha11} 
P. Zhang, \textit{Monomorphism categories, cotilting theory, and Gorenstein-projective modules}, J. Algebra \textbf{339} (2011), 181--202.

 \bibitem[Zha13]{Zha13}
    P. Zhang, \textit{Gorenstein-projective modules and symmetric recollements}, J. Algebra \textbf{388} (2013), 65--80.



	\end{thebibliography}
\end{document}